\documentclass[conference,letterpaper]{IEEEtran}
\usepackage[left=19.1mm, right=19.1mm, top=18.5mm, bottom=19.9mm]{geometry}
\IEEEoverridecommandlockouts
\usepackage{mathtools,amssymb,lipsum}

\usepackage{cuted}
\usepackage{cite}

\DeclareFontFamily{OT1}{pzc}{}
\DeclareFontShape{OT1}{pzc}{m}{it}{<-> s * [1.10] pzcmi7t}{}
\DeclareMathAlphabet{\mathpzc}{OT1}{pzc}{m}{it}
\usepackage{amsmath,amssymb,amsfonts,amsthm}
\newtheoremstyle{assumptionstyle}
   {0pt}
   {0pt}
   {}
   {11pt}
   {\em}
   {\em}
   {5pt}
   {}
\theoremstyle{assumptionstyle}

\usepackage{algorithm}
\usepackage{algorithmic}
\usepackage{graphicx}
\usepackage{textcomp}
\usepackage{xcolor}
\def\BibTeX{{\rm B\kern-.05em{\sc i\kern-.025em b}\kern-.08em
    T\kern-.1667em\lower.7ex\hbox{E}\kern-.125emX}}
\usepackage{array,multirow,graphicx} 
\usepackage{tikz}   
\def\checkmark{\tikz\fill[scale=0.3](0,.35) -- (.25,0) -- (1,.7) -- (.25,.15) -- cycle;} 

\usepackage{caption}
\usepackage{blindtext}
\usepackage{array}
\newcolumntype{P}[1]{>{\centering\arraybackslash}p{#1}}
\newcolumntype{M}[1]{>{\centering\arraybackslash}m{#1}}

\usepackage{todonotes}
\usepackage{comment}
\newcommand{\myequation}{\begin{equation}}
\newcommand{\myendequation}{\end{equation}}
\let\[\myequation
\let\]\myendequation
\allowdisplaybreaks
\usepackage{float}
\usepackage{longtable}

\begin{document}
\title{\vspace*{+6.mm} A Dynamic Sustainable Competitive Petroleum Supply Chain Model for Various Stakeholders with Shared Facilities 
\thanks{}
}

\author{Nazanin Moradinasab$^{1}$, Hassan Jafarzadeh$^{2}$, M. R. Amin-Naseri$^{3}$ and Cody H. Fleming $^{4}$
\thanks{$^{1}$ Ph.D. Industrial and Systems Engineering, Tarbiat Modares University, Tehran, Iran
        }%
\thanks{$^{2}$ Graduate Research Assistant, Engineering Systems and Environment, University of Virginia, Charlottesville, VA
        }%
\thanks{$^{3}$ Associate Professor, Industrial and Systems Engineering, Tarbiat Modares University, Tehran, Iran 
        }%
\thanks{$^{4}$Assistant Professor, Engineering Systems and Environment, University of Virginia, Charlottesville, VA 
        }%
}

\maketitle
\thispagestyle{plain} 
\pagestyle{plain} 

\begin{abstract}
Petroleum industry is the world's biggest energy source, and its associated industries such as production, distribution, refining and retail are considered as the largest ones in the world. Having the increasing price and government’s job creation and international environmental policies, the petroleum companies try to maximize the number of created job, and their profit and minimize the air pollution simultaneously. To meet these objectives, an effective detailed and precise planning is needed. On the other hand, the dynamic environment and the presence of various stakeholders add to the complexity of planning and design of petroleum supply chain. Therefore, the multi- period, multi-objective, multi-level and multi-product dynamic sustainable competitive petroleum supply chain (DSCPSC) model taking into consideration the various stakeholders have been proposed in this paper. The proposed model is an MILP model and GAMS 24.1.2 software has been used to run it for a part of real petroleum supply chain data. Sensitivity analysis was then performed to determine the sensitivity of the results to the variation of the coefficients in objective function. Sensitivity analysis reveals that the highest variations of the objective function were observed with respect to the variable costs, facility installation costs and pipeline transportation costs.

\end{abstract}

\begin{IEEEkeywords}
Competition, petroleum supply chain, sustainability, Mix integer linear problem.
\end{IEEEkeywords}

\section{Introduction}\label{sec:intro}

Petroleum supply chain (PSC) is a large network of facilities, infrastructures and complex processes from crude oil extraction to the delivery of the final product to the final customer. This network consists of downstream, midstream and upstream sectors in a wide geographical level in which there are many PSC network design decisions in strategic, tactical and operational levels. Obviously, the management and planning of such network of decisions is very complicated and difficult. Literature survey shows that much work has been done on PSC, in most of which a part of network design decisions shown in Table \ref{review_table} have been studied. As observed, some of the PSC has been investigated by applying simplifying assumptions in most works thereby a general and optimal response has not been achieved. Consideration of all aspects of the problem clearly increases the complexity of the problem and the time required to solve it. Therefore, a balance must be kept between simplifying assumptions and problem complexity.
As we know, one of the inevitable factors in PSC is the environmental regulations imposed on the network whereas only Ribas (2010) has considered environmental impacts as a constraint in his model \cite{ribas2010optimization}. On the other hand, since PSC is one of the industries requiring high investment, it can consequently lead to a high volume of employment. Therefore, consideration of employment planning in PSC is of special significance for governments. Maximization of employment and minimization of pollution are therefore important and inevitable factors together with profit maximization in PSC. Taking these three objectives into account at the same time forms the sustainable supply chain based on Carter and Rogers’ definition.

According to the Carter and Rogers’ sustainability definition, “the strategic, transparent integration and achievement of an organization’s social, environmental, and economic goals is the systematic coordination of key inter-organizational business processes for improving the long term economic performance of the individual company and its supply chain” \cite{carter2008framework}. However, only Moradinasab et al. (2018)  have considered these objects simultaneously\cite{moradinasab2018competition}.

According to Table \ref{review_table}, most researches have concentrated on planning models at the tactical level. Only Khosrojerdi et al. (2012), Nasab and Amin-Naseri (2016) and Nasab et al. (2018)  have considered both tactical and strategic levels in their studies\cite{khosrojerdi2012designing,nasab2016designing,nasab2018benders}. 

Another important factor is the presence of different stakeholders in PSC. As stated, PSC is a large network of different centers managed by various stakeholders, each looking to maximize their profits. As Table \ref{review_table} shows, only Fernandes (2013) and Moradinasab et al. (2018) have imagined different stakeholders in their works while others have considered only one stakeholder in PSC. Moradinasab et al. (2018) have developed a sustainable competitive petroleum supply chain (SCPSC) model to minimize pollution while maximizing the profits and job creation. They have considered a government and private sectors as stakeholders which have non-cooperative competition to determine their refined product prices. They have modelled the integrated petroleum supply chain by taking into account two approaches of capacity expansion including new facilities and pipelines installation 

\begin{table*}[ht]
\caption{Papers classified according to supply chain, decision levels and problem types}
\label{review_table}
\centering
\begin{tabular}{|p{3cm}|p{0.2cm}|p{0.2cm}|p{0.2cm}|p{0.2cm}|p{0.2cm}|p{0.2cm}|p{0.2cm}|p{0.2cm}|p{0.2cm}|p{0.2cm}|p{0.2cm}|p{0.2cm}|p{0.2cm}|p{0.2cm}|p{0.2cm}|p{0.2cm}|p{0.2cm}|p{0.2cm}|p{0.2cm}|p{0.2cm}|}

\hline
 &  \multicolumn{3}{|c|}{Sector} &      \multicolumn{3}{|c|}{Decision level} & \multicolumn{14}{|c|}{Supply chain decisions}\\

\cline{2-21}
 \centering Author (year)   & {\rotatebox[origin=c]{90}{Upstream}} & {\rotatebox[origin=c]{90}{Midstream}} & {\rotatebox[origin=c]{90}{Downstream}}  & {\rotatebox[origin=c]{90}{Strategic}}  & {\rotatebox[origin=c]{90}{Tactical}}&
   {\rotatebox[origin=c]{90}{Operational}} & \multicolumn{2}{|p{1.2cm}|}{Location allocation} & \multicolumn{2}{|p{1.2cm}|}{Capacity expansion}  &\multicolumn{2}{|p{1.2cm}|}{Capacity reduction}  &  {\rotatebox[origin=c]{90}{Inventory}} & {\rotatebox[origin=c]{90}{Raw material Preparation }} &{\rotatebox[origin=c]{90}{Production}} & {\rotatebox[origin=c]{90}{Routing}} & {\rotatebox[origin=c]{90}{Transportation modes}} &{\rotatebox[origin=c]{90}{Multi- Stakeholders}} & {\rotatebox[origin=c]{90}{Sharing Facilities}} & {\rotatebox[origin=c]{90}{Multi- Objective}} \\
\cline{8-13}
 
 &  &  &   &   & & & {\rotatebox[origin=c]{90}{Facility}} & {\rotatebox[origin=c]{90}{Route}}  & {\rotatebox[origin=c]{90}{Facility}} & {\rotatebox[origin=c]{90}{Route}}& {\rotatebox[origin=c]{90}{Facility}} & {\rotatebox[origin=c]{90}{Route}} &  &  &   &   & & & & \\
  \hline
  Sear 1993\cite{sear1993logistics
  } & & \checkmark& \checkmark& & & & \checkmark & & & & & & & & &\checkmark & \checkmark & & & \\
  \hline
  Escudero, Quintana et al. 1999\cite{escudero1999coro} & & \multirow{2}{*}{\centering\checkmark}& \multirow{2}{*}{\centering\checkmark}& & \multirow{2}{*}{\centering\checkmark} & &  & & & & & & & \multirow{2}{*}{\centering\checkmark}& & &  & & & \\
  \hline
    Dempster, Pedron et al. 2000 \cite{dempster2000planning}& & \multirow{2}{*}{\centering\checkmark} & \multirow{2}{*}{\centering\checkmark}& & \multirow{2}{*}{\centering\checkmark} & &  & & & & & & & & & &  & & & \\
  \hline

   Pinto, Joly et al. 2000 \cite{pinto2000planning}& & \checkmark& & &  &\checkmark &  & & & & & & \checkmark& & \checkmark& &  & & & \\
  \hline
  
    Neiro and Pinto 2006 \cite{neiro2006lagrangean}&\checkmark & \checkmark& \checkmark& & \checkmark & \checkmark&  & & & & & &\checkmark &\checkmark & \checkmark & & \checkmark & & & \\
  \hline

    Mendez, Grossmann et al. 2006\cite{mendez2006simultaneous
    } & &\multirow{2}{*}{\centering\checkmark} & & &  & \multirow{2}{*}{\centering\checkmark}&  & & & & & &\multirow{2}{*}{\centering\checkmark} & & \multirow{2}{*}{\centering\checkmark} & &  & & & \\
  \hline
  
      Pongsakdi, Rangsunvigit et al. 2006\cite{pongsakdi2006financial
      }&\multirow{2}{*}{\centering\checkmark} &\multirow{2}{*}{\centering\checkmark} & & & \multirow{2}{*}{\centering\checkmark} & &  & & & & & &\multirow{2}{*}{\centering\checkmark} & & \multirow{2}{*}{\centering\checkmark} & &  & & & \\
  \hline
  
        Al-Othman, Lababidi et al. 2008 \cite{al2008supply} &\multirow{2}{*}{\centering\checkmark} &\multirow{2}{*}{\centering\checkmark} & \multirow{2}{*}{\centering\checkmark}& & \multirow{2}{*}{\centering\checkmark} & &  & & & & & & &\multirow{2}{*}{\centering\checkmark} & \multirow{2}{*}{\centering\checkmark} & &  & & & \\
  \hline
  
          Kuo and Chang 2008\cite{kuo2008application} &\checkmark &\checkmark & \checkmark& &  &\checkmark &  & & & & & &\checkmark &\checkmark & \checkmark & &  & & & \\
  \hline
  
            Pitty, Li et al. 2008 \cite{pitty2008decision} &\checkmark & \checkmark &\checkmark &  &\checkmark&\checkmark&& & & & & &\checkmark &\checkmark & \checkmark & &  & & & \\
  \hline
  
             Kim, Yun et al. 2008\cite{kim2008integrated}& & \checkmark &\checkmark &  &\checkmark&\checkmark&\checkmark& & & & & & & &  & & \checkmark & & & \\
  \hline
  
               Al-Qahtani and Elkamel 2008\cite{al2008multisite} & & \multirow{2}{*}{\centering\checkmark} &\multirow{2}{*}{\centering\checkmark} &  \multirow{2}{*}{\centering\checkmark}&&&& &\multirow{2}{*}{\centering\checkmark} & & & & &\multirow{2}{*}{\centering\checkmark} &  & \multirow{2}{*}{\centering\checkmark}&  & & & \\
  \hline
  
                 MirHassani 2008 
                \cite{mirhassani2008operational}& &  &\checkmark &  &\checkmark&&& & & & & & \checkmark& &  \checkmark& &  & & & \\
  \hline
                   Guyonnet, Grant et al. 2008\cite{guyonnet2008integrated}&\multirow{2}{*}{\centering\checkmark} &  \multirow{2}{*}{\centering\checkmark} & \multirow{2}{*     }{\centering\checkmark} &  &&&& & & & & & \multirow{2}{*}{\centering\checkmark}& &  \multirow{2}{*}{\centering\checkmark}& &  & & & \\
  \hline
  
             Rocha, Grossmann et al. 2009\cite{rocha2009petroleum} &\multirow{2}{*}{\centering\checkmark} & \multirow{2}{*}{\centering\checkmark} &  &  &&\multirow{2}{*}{\centering\checkmark}&& & & & & & \multirow{2}{*}{\centering\checkmark}& &  \multirow{2}{*}{\centering\checkmark}& &  & & & \\
  \hline
  
              Al-Qahtani and Elkamel 2011\cite{al2011planning}& &  \multirow{2}{*}{\centering\checkmark} & \multirow{2}{*}{\centering\checkmark} & \multirow{2}{*}{\centering\checkmark} &&&& & & & & & &\multirow{2}{*}{\centering\checkmark} &  & &  & & & \\
  \hline
  
             Ribas, Hamacher et al. 2010\cite{ribas2010optimization}& &  \multirow{2}{*}{\centering\checkmark} & \multirow{2}{*}{\centering\checkmark} &  &\multirow{2}{*}{\centering\checkmark}&&& & & & & & &\multirow{2}{*}{\centering\checkmark} &  & &\multirow{2}{*}{\centering\checkmark}  & & & \\
  \hline
  
         Gill 2011\cite{gill2011supply}& &  & \checkmark &  &\checkmark&&\checkmark& & & & & & & &  & &  & & & \\
  \hline
  
         Khosrojerdi, Hadizadeh et al. 2012\cite{khosrojerdi2012designing}& & \multirow{2}{*}{\centering\checkmark} & \multirow{2}{*}{\centering\checkmark} & \multirow{2}{*}{\centering\checkmark} &\multirow{2}{*}{\centering\checkmark}&&\multirow{2}{*}{\centering\checkmark}& & & & & & \multirow{2}{*}{\centering\checkmark}& &\multirow{2}{*}{\centering\checkmark}  & &\multirow{2}{*}{\centering\checkmark}  & &\multirow{2}{*}{\centering\checkmark} & \\
  \hline
  
           Chen, Grossmann et al. 2012\cite{chen2012comparative}& & \multirow{2}{*}{\centering\checkmark} &  &  &\multirow{2}{*}{\centering\checkmark}&&& & & & & & \multirow{2}{*}{\centering\checkmark}& &  & &  & & & \\
           
  \hline
  
           Fernandes, Relvas et al. 2013\cite{fernandes2013strategic}& & & \multirow{2}{*}{\centering\checkmark}  & \multirow{2}{*}{\centering\checkmark}  &&&& & & & & & & &  & &  \multirow{2}{*}{\centering\checkmark}& \multirow{2}{*}{\centering\checkmark}& & \\
   \hline
  
          Nasab and Amin-Naseri 2016\cite{nasab2016designing}& \multirow{2}{*}{\centering\checkmark}&\multirow{2}{*}{\centering\checkmark} & \multirow{2}{*}{\centering\checkmark}  & \multirow{2}{*}{\centering\checkmark}  &\multirow{2}{*}{\centering\checkmark}&&\multirow{2}{*}{\centering\checkmark}&\multirow{2}{*}{\centering\checkmark} &\multirow{2}{*}{\centering\checkmark} &\multirow{2}{*}{\centering\checkmark} & & &\multirow{2}{*}{\centering\checkmark} & & \multirow{2}{*}{\centering\checkmark} & \multirow{2}{*}{\centering\checkmark}&  \multirow{2}{*}{\centering\checkmark}& & & \\
 
   \hline

          Moradinasab et al. 2018\cite{moradinasab2018competition}& &\multirow{2}{*}{\centering\checkmark} & \multirow{2}{*}{\centering\checkmark}  & \multirow{2}{*}{\centering\checkmark}  &\multirow{2}{*}{\centering\checkmark}&&\multirow{2}{*}{\centering\checkmark}&\multirow{2}{*}{\centering\checkmark} &\multirow{2}{*}{\centering\checkmark} &\multirow{2}{*}{\centering\checkmark} & & &\multirow{2}{*}{\centering\checkmark} & & \multirow{2}{*}{\centering\checkmark} & \multirow{2}{*}{\centering\checkmark}&  \multirow{2}{*}{\centering\checkmark}&\multirow{2}{*}{\centering\checkmark} & &\multirow{2}{*}{\centering\checkmark} \\
 
   \hline  
           Nasab et al. 2018 \cite{nasab2018benders}& \multirow{2}{*}{\centering\checkmark}& \multirow{2}{*}{\centering\checkmark}& \multirow{2}{*}{\centering\checkmark}  & \multirow{2}{*}{\centering\checkmark}  &\multirow{2}{*}{\centering\checkmark}&&\multirow{2}{*}{\centering\checkmark}&\multirow{2}{*}{\centering\checkmark} &\multirow{2}{*}{\centering\checkmark} &\multirow{2}{*}{\centering\checkmark} & & &\multirow{2}{*}{\centering\checkmark} & &  \multirow{2}{*}{\centering\checkmark}&\multirow{2}{*}{\centering\checkmark} & \multirow{2}{*}{\centering\checkmark} \multirow{2}{*}& & & \\

  \hline
\end{tabular}
\end{table*}

and capacity expansion of the present facilities and pipelines. Integrated and simultaneous consideration of these two approaches is important in PSC\cite{moradinasab2018competition}. According to Table \ref{review_table}, Al-Qahtani and Elkamel (2008) have studied only the capacity expansion of the facilities in PSC \cite{al2008supply}. In addition, only 7 references have considered location and allocation of the facilities in the midstream sector (Table \ref{review_table}). As observed in the Table, the feature of reduced capacity or close down of facilities has not been considered in any of the researches. However, one of the important features in PSC is the possibility of close down or capacity reduction during different time periods. The major and important reason for this consideration is the high fixed and variable costs of existing facilities compared with new ones. This assumption makes the capacity reduction and close down of some of the existing facilities possible by increasing the PSC capacity via constructing new facilities and expansion of existing facilities during different time periods. Simultaneous consideration of all three approaches of new facilities construction, existing facilities expansion and capacity reduction and/or close down of the existing facilities is referred to as the dynamic design of supply chain \cite{arabani2012facility}. Athough the dynamic design of supply chain is important in the petroleum supply chain supply chain because of the high fixed and variable costs for facilities, it is not investigated in any of the sudies which have done so far.
 
 In this paper, multi- period, multi objective, multi-level and multi-product dynamic sustainable competitive petroleum supply chain (DSCPSC) in the presence of various stakeholders with shared facilities in the three levels of refineries, distribution centers and customers is proposed. Being strategic, the petroleum upstram sector including crude oil fields and storage teminals have not been consieded in this study. Since these sectors are critical and strategic for governments, they are kept under governement control and not privatized (i.e. they are ignored in the competitive petroleum supply chain). The decision variables in this model including location, new facilities and pipelines installation, capacity expansion in the existing facilities and pipelines, close down of facilities, inventory, production, imports, exports, routing and selection of transportation modes in strategic and tactical levels. In the model presented at the next sections the routing part has been defined as the shortest path between defined centers, while to have a more detailed model shortest traveling time and also the congestion of the routes can be added to it as well \cite{benzaman2019discrete}. 
 
In this work all the three aspects including installation, closing down and capacity expansion of the facilities and pipelines are studied at the same time which makes the presented model dynamic. Another aspect of the dynamism in the proposed model is associated with the distribution centers. Distribution centers are included in the shared resources of stakeholders, that means different stakeholders can use the same distribution center within the same time period, but the percentage of using a distribution center by stockholders can change within different time periods. Each stakeholder pays the costs of using a distribution center based on the percentage of using the center. In addition to optimal use of the capacity of a distribution center, one of the advantages of this approach is the minimization of each stakeholder’s expenses.

The optimization model generated for this problem is in MILP form, because it has shown its strength in providing a convenient way to capture most aspects of the different problems. That is why it is utilized in different modeling research fields such as path planning \cite{jafarzadeh2017enhanced}, scheduling \cite{moradinasab2013no}, logistics \cite{jafarzadeh2017genetic}, dynamic control problems \cite{earl2005iterative}, chemical processing \cite{huang2018water}, pricing and inventory Control \cite{rasouli2014joint} and etc. Another advantage of MILP models is the existence of powerful solvers in the market such as CPLEX. A downside of MILP is the complexity of the model. Using integer variables makes the model NP-hard and increases the size of solution space exponentially. A technique to overcome this weakness is applying heuristic algorithms along with traditional MILP solving methods. For instance, rapid path planning algorithms such as \cite{jafarzadeh2014new, jafarzadeh2018exact} can be used when a part of the modeling are dealing with finding a shortest path between different facilities.

In this paper, the problem definition is given in the next section. Following the statement of problem assumptions, modeling is dealt with in section 4. Fuzzy planning for solving multi-objective models is presented in section 5. Finally, section 6 gives the numerical results obtained by executing the model and conclusions and suggestions for future researches are provided in section 7.



\section{Problem Definition}\label{sec:prdif}

In this work, petroleum supply chain is studied in three levels including refineries, distribution centers, and customer centers in presence of different stakeholders. Having various stakeholders in a network leads to competition between them. In this network, the demand of each customer domain can be met by different stakeholders which in this study are different supply chains.  For optimal use of distribution centers, each of them is assumed to be owned by several stakeholders and each stakeholder paying the fixed and variable costs of the distribution center according to the portion of center capacity usage. 

Profit maximization is one of the most common objective functions according to the literature survey. The environmental regulations and the significance of the increasing job opportunities make the minimization of the pollution, maximization of profit and job creation important for governments. These have not been paid enough attention in the literature, therefore this study emphasizes more on these objectives. The supply chain in which the objectives of minimizing the pollution, maximizing the job opportunities and maximizing the overall profit are modeled simultaneously, is called sustainable supply chain. Having different stakeholders, the DSCPSC model is used to study the petroleum supply chain in this work. The DSCPSC model optimizes the supply chain decisions for each stakeholder simultaneously, i.e. optimal locations for installing facilities (refineries, distribution centers and pipelines), capacities expansions, the selection between the two installation and capacity expansion policies (refineries, distribution centers and pipelines), close down of the existing facilities (refineries, distribution centers and pipelines), inventory, pipeline route, selection of transportation modes, the amounts of crude oil and refined product and the amounts of import and exports.

In this work, different mode of transportation is available for stakeholders. They transport refined products such as gasoline, gasoil and kerosene from refineries to distribution centers by one of the transportation modes and ultimately sells it to the final customers.

\section{Problem Assumptions}\label{sec:prasum}

This section explains the problem assumption in DSCPSC mathematical model. The following assumptions help to convert a real problem to the proposed mathematical model:

\begin{enumerate}
    \item The supply chain consists of the three levels of refineries, distribution centers, and final customers.
    \item Customer demands are definite and clear.
    \item The supply chain includes several levels and products and different transportation modes.
    \item The location of the existing refineries, distribution centers and existing pipelines is fixed and may not be changed.
    \item Capacity expansion may be performed in existing refineries, distribution centers and pipelines and not in new facilities and pipelines.
    \item	Closing down the facilities may occur in existing refineries, distribution centers and existing pipelines and not in new facilities and pipelines.
    \item The stream of refined products is categorized into three parts in this supply chain. These categories and the corresponding modes of transportation are as follows:
    \begin{itemize}
     \item Transportation between refineries and distribution centers, which are carried out by roads, rail and pipelines.
     \item Transportation between distribution centers, which are carried out by roads, rail and pipelines.
     \item Transportation between distribution centers and final customers, which are carried out by roads and rail.
   \end{itemize}
    \subitem 
•	
    \item New pipelines may be constructed in the first and second transportation categories.
    \item The number and type of storage tanks for each refined product in each distribution center are determined according to the input stream into the distribution center.
    \item The distribution centers may either export or import the refined products.
    11. The costs of capacity expansion, construction and closing down an existing facility          and the fixed and variable costs for each type of facility are known given its capacity.
  \item  The inventory holding costs for each of the facilities in each time period have been considered in this supply chain.
   \item A part of the inventory is kept as the confidence inventory in each facility to prevent shortages.
   \item  All possible levels of capacity for new refineries, distribution centers and pipelines are known.
  \item  All possible levels of capacity expansion for existing refineries, distribution centers and crude oil pipelines are known. 
   \item  In addition to economical profit function, two other functions including minimization of pollution and maximization of job creation have also been considered.

\end{enumerate}

\section{Problem Modeling}\label{sec:prmodle}

The formulation and mathematical DSCPSC model, an MILP model, is presented in this part. The definitions of the parameters, sets and variables used in the model are given in Appendix A.

\subsection{Objective function}

DSCPSC model is modeled as a multi-objective linear programming model with three objective functions. The first objective function is profit maximization and the second objective function is the minimizing the amount of pollutions resulting from facilities and transportation modes. The third objective function is maximizing the number of created job. The objective functions are formulated as follows:

\subsubsection{Total profit maximization} 
According to Equation \ref{Multi1}, each stakeholder’s profit is obtained from the difference between incomes and expenses. The incomes for each of the stakeholders, obtained by Equation \ref{Multi2}, includes those from the sale and export of the products.

  \begin{flalign} \label{Multi1}
     Maximize \; P_{e}=RF_{e}-CF_{e} && 
  \end{flalign}

  \begin{flalign} \label{Multi2}
     RF_{e}&=\sum_{t\in T}\sum_{p\in P}\sum_{v\in V} \sum_{l\in (L_{e}\cup \acute{L_{e}})}\sum_{m\in M} qlm_{p,l,m,v}^{t,e}pr_{}^{e,p,t}\nonumber\\ 
     &+\sum_{t\in T}\sum_{p\in P} \sum_{l\in (L_{e}\cup \acute{L_{e}})}Ep_{p,l}^{t,e}ERPP_{p}^{t,e} 
  \end{flalign}

According to Equation \ref{Multi3}, which shows each of the stakeholder’s expenses, the first and second relations represent the installation cost of new refineries and distribution centers, respectively. Capacity expansion costs of existing refineries and distribution centers are given by relations 3 and 4. Closing down cost of existing refineries are shown by relation 5 and relations 6-8 show installation costs of new pipelines. Capacity expansion cost of existing pipelines are calculated using relations 9-10. Installation cost of new storage tanks in distribution centers are obtained by relation 11.
The facilities constant and variable costs are given by relations 12-17. Relations 18-19 are used to obtain inventory costs of new and existing refineries and distribution centers. The costs of refined products transported between refineries and distribution centers to the final customers (i.e. is done by one of the transportation modes) are determined using relations 20-25. Product importation cost are calculated using relation 26. Finally, employment costs for refineries and distribution centers are obtained using relation 27.


\begin{strip}
  \begin{flalign} \label{Multi3}
     C_{e\;\forall e\,\in E} &=\sum_{t\in T} \sum_{k\in  \acute{K_{e}}}\sum_{ek\in EK} xcostk_{k}^{ek,t}xk_{k}^{ek,t,e}
     +\sum_{t\in T}\sum_{el\in EL} \sum_{l\in  \acute{L_{e}}} xcostl_{l}^{el,t}xl_{l}^{el,t,e}+\sum_{t\in T}\sum_{uk\in UK} \sum_{k\in K_{e}} ucostk_{k}^{uk,t}\tau k_{k}^{uk,t,e}\nonumber\\
     &+\sum_{t\in T}\sum_{ul\in UL}\sum_{p\in P} \sum_{l\in L_{e}} ucostl_{l}^{ul,p,t}\tau l_{p,l}^{ul,t,e}+\sum_{t\in T}\sum_{uk\in UK} \sum_{k\in K_{e}} clcostk_{k}^{t}\psi k_{k}^{t,e}+\sum_{t\in T}\sum_{lv\in LV} \sum_{k\in (K_{e}\cup\acute{K_{e}})}\sum_{l\in (L_{e}\cup\acute{L_{e}})}\nonumber\\
     &\sum_{v=3} rcostkl_{k,l,v}^{lv,rv,t}r_{k,l,v}^{lv,rv,t,e}+\sum_{t\in T} \sum_{lp\in L_{e}}\sum_{l\in (L_{e}\cup\acute{L_{e}})}\sum_{lv\in LV}\sum_{v=3}\sum_{rv\in RV} rcostlpl_{lp,l,p}^{lv,rv,t}rlpl_{lp,l,v}^{lv,rv,t,e}+\sum_{t\in T} \sum_{lp\in \acute{L_{e}}}\sum_{l\in \acute{L_{e}}}\nonumber\\
     &\sum_{lv\in LV}\sum_{v=3}\sum_{rv\in RV} rcostlpl_{lp,l,p}^{lv,rv,t}rlpl_{lp,l,v}^{lv,rv,t,e}+\sum_{t\in T} \sum_{l\in L_{e}}\sum_{k\in K_{e}}\sum_{ev\in EV} ycostkl_{k,l}^{ev,t}ykl_{k,l}^{ev,t,e}+\sum_{t\in T} \sum_{l\in L_{e}}\sum_{lp\in L_{e}}\sum_{ev\in EV} \nonumber\\
     & ycostlpl_{lp,l}^{ev,t}ylpl_{lp,l}^{ev,t,e}+\sum_{t\in T} \sum_{l\in \acute{L_{e}}}\sum_{p\in P}\sum_{ez\in EZ} ncost_{l}^{ez,t}n_{p,l}^{ez,t,e}+\sum_{t\in T} \sum_{l\in L_{e}}\sum_{p\in P} Fcostl_{l}^{t}(\xi_{l}^{t,e}icl_{p,l}^{}+\sum_{ul\in UL}capl_{l}^{ul}\nonumber\\
     &\tau l_{p,l}^{ul,t,e})+\sum_{t\in T} \sum_{k\in K_{e}} Fcostk_{k}^{t}(ick_{k}^{}(1-\sum_{\acute{t}\leqslant t\in T}\psi k_{k}^{\acute{t},e})+\sum_{uk\in UK}capk_{k}^{uk}\tau k_{k}^{uk,t,e})
    +\sum_{t\in T} \sum_{k\in \acute{K_{e}}}\sum_{ek\in EK} Fcostk_{k}^{t}\nonumber\\
    &xk_{k}^{ek,t,e}c_{k}^{ek}+\sum_{t\in T} \sum_{l\in \acute{L_{e}}}\sum_{ez\in EZ}\sum_{p\in P}  Fcostl_{l}^{t}Nct_{t}^{ez}n_{p,l}^{ez,l,e}+\sum_{t\in T} \sum_{l\in \acute{L_{e}}}\sum_{v\in V}\sum_{p\in P} (\sum_{k\in (K_{e}\cup\acute{K_{e}})} qkl_{p,k,l,v}^{t,e}+\nonumber\\
    &\sum_{lp\in \L_{e}}qlpl_{p,lp,l,v}^{t,e}-\sum_{lp\in (L_{e}\cup\acute{L_{e}})} qlpl_{p,l,lp,v}^{t,e}) pcostl_{p,l}^{t}+\sum_{t\in T} \sum_{l\in \acute{L_{e}}}\sum_{v\in V}\sum_{p\in P} (\sum_{k\in (K_{e}\cup\acute{K_{e}})}
    qkl_{p,k,l,v}^{t,e}+\sum_{lp\in (L_{e}\cup\acute{L_{e}})}\nonumber\\
    & qlpl_{p,lp,l,v}^{t,e}-\sum_{lp\in \acute{L_{e}}} qlpl_{p,l,lp,v}^{t,e}) pcostl_{p,l}^{t}+\sum_{t\in T} \sum_{k\in (K_{e}\cup\acute{K_{e}})} OP_{}^{t} hcostk_{k}^{}vk_{k}^{t,e}+ \sum_{p\in P}\sum_{l\in (L_{e}\cup\acute{L_{e}},e\in E_{PS})}
    pr_{}^{e,p,t}\nonumber\\
    & hcostl_{p,l}^{}vl_{p,l}^{t,e}+\sum_{t\in T} \sum_{k\in (K_{e}\cup\acute{K_{e}})}\sum_{v\in V}\sum_{p\in P}\sum_{l\in (L_{e}\cup\acute{L_{e}})}
    qkl_{p,k,l,v}^{t,e} qcostkl_{k,l}^{t}+\sum_{t\in T}\sum_{p\in P}\sum_{l\in (L_{e}\cup\acute{L_{e}})}\sum_{lp\in \acute{L_{e}}}qlpl_{p,lp,l,v}^{t,e} \nonumber\\
    &qcostlpl_{lp,l}^{t}
    +\sum_{t\in T}\sum_{p\in P}\sum_{l\in \acute{L_{e}}}\sum_{lp\in \acute{L_{e}}}qlpl_{p,lp,l,v}^{t,e} qcostlpl_{lp,l}^{t}
    +\sum_{t\in T}\sum_{v\in V}\sum_{l\in (L_{e}\cup\acute{L_{e}})}\sum_{k\in (K_{e}\cup\acute{K_{e}})}\sum_{lcv\in LCV} ncostkl_{k,l,v}^{lcv,t} \nonumber\\
    &nkl_{p,k,l,v}^{lcv,t,e}
    +\sum_{t\in T}\sum_{v\in V}\sum_{l\in (L_{e}\cup\acute{L_{e}})}\sum_{m\in M}\sum_{p\in P}\sum_{lcv\in LCV}ncostlm_{l,m,v}^{lcv,t} nlm_{p,l,m,v}^{lcv,t,e}
    +\sum_{t\in T}\sum_{v\in V}\sum_{l\in (L_{e}\cup\acute{L_{e}})}\sum_{lp\in (L_{e}\cup\acute{L_{e}})}\nonumber\\
    &\sum_{lcv\in LCV}\sum_{p\in P}ncostlpl_{l,lp,v}^{lcv,t} nlpl_{p,lp,l,v}^{lcv,t,e}
    +\sum_{t\in T}\sum_{p\in P}\sum_{l\in (L_{e}\cup\acute{L_{e}})}i_{p,l}^{t,e}icost_{p}^{t}+\sum_{t\in T}\sum_{lev\in LEV}WCost_{lev}^{t}[\sum_{k\in \acute{K_{e}}}\sum_{\acute{en}\in EN}\nonumber\\
    &\sum_{en\in EN}HENK_{en,k}^{\acute{en},lev,t}
    +\sum_{k\in K_{e}}\sum_{\acute{en}\in EN}\sum_{en\in EN}HEEK_{en,k}^{\acute{en},lev,t}+\sum_{l\in \acute{L_{e}}}\sum_{\acute{en}\in EN}\sum_{en\in EN}HENL_{en,l}^{\acute{en},lev,t}+\sum_{l\in L_{e}}\nonumber\\
    &\sum_{\acute{en}\in EN}\sum_{en\in EN}HEEL_{en,l}^{\acute{en},lev,t}]
    \nonumber\\
    &
    \end{flalign}
\end{strip}

\subsubsection{Minimizing the pollution created by facilities and transportation modes}

In Equation \ref{Multi4}, the first and second terms show the pollution by new refineries and distribution centers, respectively. Parameter $\lambda E_{en}^{}$ is an indication of the sensitivity of an area to the pollution. The amount of pollution by capacity expansion of existing refineries and distribution centers are calculated by the fourth and fifth terms in Equation 4, respectively. Finally, the pollution by transportation modes are obtained by the last two terms in Equation \ref{Multi4}.

    \begin{flalign} \label{Multi4}
     Minimize \; Pul_{e}&=
     \sum_{t\in T}\sum_{k\in \acute{K_{e}}}\sum_{en\in EN}\sum_{ek\in EK}\lambda E_{en}^{}Nk_{en}^{k} c_{k}^{ek}\nonumber\\ 
     &xk_{k}^{ek,t,e}Pulk+\sum_{t\in T}\sum_{l\in \acute{L_{e}}} \sum_{en\in EN}\sum_{ez\in EZ}\sum_{p\in P}
     \nonumber\\ 
     &\lambda E_{en}^{}Nl_{en}^{l} c_{l}^{ez}n_{p,l}^{ez,t}Pull+\sum_{t\in T}
     \sum_{en\in EN}\sum_{uk\in UK}
     \nonumber\\ 
     &\sum_{k\in K_{e}}
     (\lambda E_{en}^{}
    Nek_{en}^{k} c_{k}^{uk}\tau_{k}^{uk,t,e}Pulk)Per+\nonumber\\ 
     &\sum_{t\in T}
     \sum_{en\in EN}\sum_{p\in P}\sum_{ul\in UL}\sum_{l\in L_{e}}
     (\lambda E_{en}^{}Nel_{en}^{l}c_{l}^{ul}
    \nonumber\\ 
     &\tau_{p,l}^{ul,t,e}Pull)Per+\sum_{t\in T}\sum_{l\in (L_{e}\cup\acute{L_{e}})}\sum_{m\in M}
     \nonumber\\ 
     &\sum_{lcv\in LCV}\sum_{v\in V}Pulv_{lcv}^{v}nlm_{p,l,m,v}^{lcv,t,e}dis_{m}^{l}
     \nonumber\\ 
     &
  \end{flalign}
  
\subsubsection{Maximization of job creation}  

Based on Equation \ref{Multi5}, created job opportunities by installing new refineries and distribution centers are calculated by the first two terms and those by capacity expansion of the existing refineries and distribution centers are obtained using the second two terms.

      \begin{flalign} \label{Multi5}
     Maximize& \; S_{e}=\nonumber \\
     &\sum_{lev\in LEV}\sum_{k\in \acute{K_{e}}}\sum_{en\in EN}\sum_{\acute{en}\in EN}Nk_{en}^{k}W_{en}^{\acute{en}} \nonumber\\ 
     &HENK_{en,k}^{\acute{en},lev}+\sum_{lev\in LEV}\sum_{k\in {K_{e}}}\sum_{en\in EN}\sum_{\acute{en}\in EN}
     \nonumber\\ 
     &Nek_{en}^{k}W_{en}^{\acute{en}}HEEK_{en,k}^{\acute{en},lev}+\sum_{lev\in LEV}\sum_{l\in \acute{L_{e}}}
     \nonumber\\ 
     &\sum_{en\in EN}\sum_{\acute{en}\in EN}Nl_{en}^{l}W_{en}^{\acute{en}}HENl_{en,l}^{\acute{en},lev}+
     \nonumber\\ 
     &
     \sum_{lev\in LEV}\sum_{l\in {L_{e}}}\sum_{en\in EN}\sum_{\acute{en}\in EN}Nel_{en}^{l}W_{en}^{\acute{en}}\nonumber\\ 
     &HEEl_{en,l}^{\acute{en},lev}
     \nonumber\\ 
     &
  \end{flalign}
  \subsection{Model constraints}
  Several constraints have been considered in this model, as follows:
  \subsubsection{Network design constraints }

  \begin{flalign} \label{cons6}
     \sum_{t\in T} \sum_{el\in EL}xl_{l}^{el,t,e}\leq1\;\;\;\;\;\forall e \in E ,l \in \acute{L_{e}^{}}&& 
  \end{flalign}
  
 \begin{flalign} \label{cons7}
     \sum_{t\in T} \sum_{ek\in EK}xk_{k}^{ek,t,e}\leq1\;\;\;\;\;\forall e \in E ,k \in \acute{K_{e}^{}}&& 
  \end{flalign}
  
 \begin{flalign} \label{cons8}
     \sum_{t\in T} z_{l}^{ez,t,e}\leq1\;\;\;\;\;\;\;\forall ez \in EZ ,l \in \acute{L_{e}^{}},p \in P&& 
  \end{flalign}
  
  \begin{flalign} \label{cons9}
     M\times z_{l}^{ez,t,e}\geq n_{p,l}^{ez,t,e}\;\;\;\forall ez \in EZ ,l \in \acute{L_{e}^{}},p \in P,e \in E&& 
  \end{flalign}
  
  \begin{flalign} \label{cons10}
     \sum_{t\in T}\sum_{uk\in UK}\tau k_{k}^{uk,t,e}\leq1\;\;\;\;\;\;\; \forall k \in K_{e}^{} ,e \in E&& 
  \end{flalign}
  
  \begin{flalign} \label{cons11}
     \sum_{t\in T}\sum_{ul\in UL}\sum_{e\in E_{l}^{}}\tau l_{p,l}^{ul,t,e}\leq1\;\;\;\;\;\;\; \forall l \in L_{e}^{} ,p \in P, e \in E_{l}^{}&& 
  \end{flalign}
  
  \begin{flalign} \label{cons12}
     \sum_{t\in T}\sum_{ev\in EV}ykl_{k,l}^{ev,t}\leq1\;\;\;\;\;\;\; \forall l \in L_{e}^{} ,k \in K_{e}^{}, e \in E&& 
  \end{flalign}

  \begin{flalign} \label{cons13}
     \sum_{t\in T}\sum_{ev\in EV}ylpl_{lp,l}^{ev,t}\leq1\;\;\;\;\;\;\; \forall lp\&l \in L_{e}^{} , e \in E_{l}^{}&& 
  \end{flalign}  
  
    \begin{flalign} \label{cons14}
     \sum_{t\in T}\psi k_{k}^{t,e}\leq1\;\;\;\;\;\;\;\;\;\;\;\; \forall k \in K_{e}^{} , e \in E&& 
  \end{flalign}  
  
  \begin{flalign} \label{cons15}
    \xi_{l}^{e}\leq1\;\;\;\;\;\;\;\;\;\;\;\; \forall l \in L_{e}^{} , e \in E,t \in T&& 
  \end{flalign}

  \begin{flalign} \label{cons16}
    \sum_{e\in E}\xi_{l}^{e}\leq1\;\;\;\;\;\;\;\;\;\;\;\; \forall l \in L_{e}^{} , t \in T&& 
  \end{flalign}

  \begin{flalign} \label{cons17}
    \sum_{t\in T}\sum_{uk\in UK}\tau k_{k}^{uk,t,e}\leq M(1-\sum_{t\in T}\psi k_{k}^{t,e})\;\;\; \forall k \in K_{e}^{} , e \in E&& 
  \end{flalign}

  \begin{flalign} \label{cons18}
    \sum_{p\in P}\sum_{v\in V}\sum_{l\in (L_{e}\cup\acute{L_{e}})}qkl_{p,k,l,v}^{t,e}(1-\mu_{p}^{})+\nu k_{k}^{t-1,e}\leq ick_{k}^{}
    \nonumber\\
    (1-\sum_{\acute{t}\leq t\in T}
    \psi k_{k}^{t,e})+\sum_{\acute{t}\leq t\in T}\sum_{uk\in UK}capk_{k}^{uk}\tau k_{k}^{uk,\acute{t},e}\;\;\;\;\;\;\;\;
    \nonumber\\
    \;\;\; \forall k \in K_{e}^{} , e \in E_{k}^{},t\in T&& 
  \end{flalign}

  \begin{flalign} \label{cons19}
    \sum_{e\in E_{l}^{}}\sum_{v\in V}\sum_{lp\in L_{e}}qlpl_{p,lp,l,v}^{t,e}+\sum_{e\in E_{l}^{}}\sum_{v\in V}\sum_{k\in (K_{e}\cup\acute{K_{e}})}qkl_{p,k,l,v}^{t,e}\nonumber\\
    -\sum_{e\in E_{l}^{}}\sum_{v\in V}\sum_{l\in (L_{e}\cup\acute{L_{e}})}qlpl_{p,l,lp,v}^{t,e}+\sum_{e\in E_{l}^{}}i_{p,l}^{t,e}+\sum_{e\in E_{l}^{}}vl_{p,l}^{t-1,e}
    \nonumber\\
    \leq \sum_{e\in E_{l}^{}}icl_{p,l}^{}\xi_{l}^{t,e}
    +\sum_{e\in E_{l}^{}}\sum_{\acute{t}\leq t\in T}\sum_{ul\in UL}capl_{l}^{ul}\tau l_{p,l}^{ul,\acute{t},e}\;\;\;\;\;\;\;\;\;\;
    \nonumber\\
    \;\;\; \forall l \in L_{e}^{} , t\in T&& 
  \end{flalign}

  \begin{flalign} \label{cons20}
    \sum_{\acute{t}\leq t\in T}\sum_{ek\in EK}Mk_{k}^{ek}Nck_{k}^{ek}xk_{k}^{ek,\acute{t},e}
    \leq \sum_{v\in V}\sum_{p\in P}\sum_{l\in (L_{e}\cup\acute{L_{e}})}
    \nonumber\\
    qkl_{p,k,l,v}^{t,e}(1\setminus \mu_{p}^{})\leq\sum_{\acute{t}\leq t\in T}\sum_{ek\in EK}Nck_{k}^{ek}xk_{k}^{ek,\acute{t},e}\;\;\;\;\;\;\;\;\;
    \nonumber\\
    \;\;\; \forall k \in \acute{K_{e}{}} , t\in T,e\in E_{k}^{}&& 
  \end{flalign}

  \begin{flalign} \label{cons21}
    \sum_{\acute{t}\leq t\in T}\sum_{ez\in EZ}Ml_{p,l}^{ez}Nct_{l}^{ez}n_{p,l}^{ez,\acute{t},e}
    \leq \sum_{v\in V}\sum_{k\in (K_{e}\cup\acute{K_{e}})}qkl_{p,k,l,v}^{t,e}
    \nonumber\\
    +i_{l}^{p,t,e}+\sum_{v\in V}\sum_{lp\in (L_{e}\cup\acute{L_{e}})}qlpl_{p,lp,l,v}^{t,e}-\sum_{v\in V}\sum_{lp\in \acute{L_{e}}}qlpl_{p,l,lp,v}^{t,e}\nonumber\\\leq\sum_{\acute{t}\leq t\in T}\sum_{ez\in EZ}Nct_{l}^{ez}n_{p,l}^{ez,\acute{t},e}\;\;\;\;\;\;\;\;\;\;\;\;\;\;\;\;\;\;\;\;\;\;\;\;\;\;\;\;\;\;\;\;\;\;\;\;\;\;\;\;\;\;\;\;\;
    \nonumber\\
    \;\;\; \forall l \in \acute{L_{e}{}} ,e\in E_{l}^{},t \in T,p \in P&& 
  \end{flalign}

  \begin{flalign} \label{cons22}
    \sum_{ez\in EZ}n_{p,l}^{ez,t,e}
    \leq M\sum_{el\in EL}x_{l}^{el,t,e}\;\;\;\;\;\;\;\;\;\;\;\;\;\;\;\;\;\;\;\;\;\;\;\;\;\;\;\;\;\;\;\;\;\;\;\;\;
    \nonumber\\
    \;\;\;\;\;\;\;\;\;\;\;\;\;\;\;\;\;\;\;\;\;\;\;\;\;\;\;\;\;\;\;\;\;\;\;\;\;\;\;\;\;\; \forall l \in \acute{L_{e}{}} ,e\in E_{l}^{},t \in T,p \in P&& 
  \end{flalign}

  \begin{flalign} \label{cons23}
    qkl_{p,k,l,v}^{t,e}\leq \sum_{lcv\in LCV}trc_{v}^{lcv}nkl_{p,k,l,v}^{lcv,t,e}TPP+\sum_{\acute{t}\leq t\in T}\sum_{lv\in LV}
    \nonumber\\
    \sum_{rv\in RV}clv_{lv}^{}rkl_{k,l,v}^{lv,rv,\acute{t},e}+Rkl_{k,l}^{}(icapkl_{k,l}^{}+\sum_{\acute{t}\leq t\in T}\;\;\;\;\;\;\;\;\;\;\;\;\nonumber\\\sum_{ev\in EV}capkl_{k,l}^{ev}ykl_{k,l}^{ev,\acute{t},e}\;\;\;\;\;\;\;\;\;\;\;\;\;\;\;\;\;\;\;\;\;\;\;\;\;\;\;\;\;\;\;\;\;\;\;\;\;\;\;\;\;\;\;\;\;\;\;\;\;\;\;\;\nonumber\\\forall l \in L_{e}{},k \in K_{e}{} ,e\in E_{l}^{},t \in T,p \in P,v \in V&& 
  \end{flalign}  
  
  \begin{flalign} \label{cons24}
    \sum_{\acute{t}\leq t\in T}\sum_{lv\in LV}\sum_{rv\in RV}\sum_{v\in V}rkl_{k,l,v}^{lv,rv,\acute{t},e}+\sum_{v\in V}rk\acute l_{k,l,v}^{t,e}
    \leq 1\;\;\;\;\;\;
    \nonumber\\
    \;\;\;\;\;\;\;\;\;\;\;\;\;\;\;\;\;\;\;\;\;\;\;\;\;\;\;\;\;\;\;\; \forall l \in L_{e}{},k \in K_{e}{}  ,e\in E,t \in T&& 
  \end{flalign}

  \begin{flalign} \label{cons25}
    \sum_{\acute{t}\leq t\in T}\sum_{lv\in LV}\sum_{rv\in RV}\sum_{v\in V}rkl_{k,l,v}^{lv,rv,\acute{t},e}+Rkl_{k,l}^{}
    \leq 1\;\;\;\;\;\;\;\;\;\;\;\;
    \nonumber\\
    \;\;\;\;\;\;\;\;\;\;\;\;\;\;\;\;\;\;\;\;\;\;\;\;\;\;\;\;\;\;\;\;\;\;\;\;\;\; \forall l \in L_{e}{},k \in K_{e}{}  ,e\in E,t \in T&& 
  \end{flalign}

  \begin{flalign} \label{cons26}
    nkl_{p,k,l,v}^{lcv,t,e}
    \leq M\sum_{v\in V}rk\acute l_{k,l,v}^{t,e}\;\;\;\;\;\;\;\;\;\;\;\;\;\;\;\;\;\;\;\;\;\;\;\;\;\;\;\;\;\;\;\;\;\;\;\;\;\;\;\;\;\;\;\;\;\;
    \nonumber\\
    \;\;\;\;\;\;\;\;\;\;\;\;\;\;\forall l \in L_{e}{},k \in K_{e}{}  ,e\in E,t \in T,p \in P,lcv \in LCV&& 
  \end{flalign}

  \begin{flalign} \label{cons27}
    \sum_{\acute{t}\leq t\in T}\sum_{ev\in EV}ykl_{k,l}^{ev,\acute{t},e}\leq Rkl_{k,l}^{}
     \;\;\;\;\;\;\;\;\;\;\;\;\;\;\;\;\;\;\;\;\;\;\;\;\;\;\;\;\;\;\;\;\;\;
    \nonumber\\
    \;\;\;\;\;\;\;\;\;\;\;\;\;\;\;\;\;\;\;\;\;\;\;\;\;\;\;\;\;\;\;\;\;\;\;\;\;\; \forall l \in L_{e}{},k \in K_{e}{}  ,e\in E,t \in T&& 
  \end{flalign}

  \begin{flalign} \label{cons28}
    rkl_{k,l,v}^{lv,rv,t,e}=0
     \;\;\;\;\;\;\;\;\;\;\;\;\;\;\;\;\;\;\;\;\;\;\;\;\;\;\;\;\;\;\;\;\;\;\;\;\;\;\;\;\;\;\;\;\;\;\;\;\;\;\;\;
    \nonumber\\
    \;\;\;\;\;\;\;\;\;\;\;\;\; \forall l \in L_{e}{},k \in K_{e}{}  ,e\in E,t \in T,p \in P,v \in V \nonumber\\
    v\neq pipeline\;\;route\;\; transportation\;\;mode&&
  \end{flalign}

  \begin{flalign} \label{cons29}
    rk\acute{l}_{k,l,v}^{t,e}=0
     \;\;\;\;\;\;\;\;\;\;\;\;\;\;\;\;\;\;\;\;\;\;\;\;\;\;\;\;\;\;\;\;\;\;\;\;\;\;\;\;\;\;\;\;\;\;\;\;\;\;\;\;
    \nonumber\\
    \;\;\;\;\;\;\;\;\;\;\;\;\; \forall l \in L_{e}{},k \in K_{e}{}  ,e\in E,t \in T,p \in P \nonumber\\
    v= pipeline\;\;route\;\; transportation\;\;mode&&
  \end{flalign}

   \begin{flalign} \label{cons30}
    qkl_{p,k,l,v}^{t,e} \leq\sum_{lcv\in LCV}trc_{v}^{lcv}nkl_{p,k,l,v}^{lcv,t,e}
    TPP+\sum_{\acute{t}\leq t\in T}\sum_{lv\in LV}
    \nonumber\\
    \sum_{rv\in RV}clv_{lv}^{}rkl_{k,l,v}^{lv,rv,\acute{t},e}\;\;\;\; \forall l \in (\acute{L_{e}{}\cup L_{e}{})},\forall k \in \acute{K_{e}{}} ,e\in E,
    \nonumber\\
    t \in T,p \in P,v \in V\;\;\; OR \;\;\;\forall l \in \acute{L_{e}{} },\forall k \in \acute{K_{e}{}} ,e\in E
    \nonumber\\
    ,t \in T,p \in P,v \in V \;\;\;\;\;\;\;\;\;\;\;\;\;\;\;\;\;\;\;\;\;&& 
  \end{flalign}

  \begin{flalign} \label{cons31}
    \sum_{\acute{t}\leq t\in T}\sum_{lv\in LV}\sum_{rv\in RV}\sum_{v\in V}rkl_{k,l,v}^{lv,rv,\acute{t},e}+\sum_{v\in V}rk\acute l_{k,l,v}^{t,e}
    \leq 1\;\;\;\;\;\;\;\;\;\;\;\;
    \nonumber\\
     \forall l \in (\acute{L_{e}{}\cup L_{e}{})},\forall k \in \acute{K_{e}{}} ,e\in E,t \in T\;\;\;\;\;
     \nonumber\\
     OR\;\;\forall l \in \acute{L_{e}{}},\forall k \in K_{e}{} ,e\in E,t \in T\;\;\;\;\;
     && 
  \end{flalign}

  \begin{flalign} \label{cons32}
    nkl_{p,k,l,v}^{lcv,t,e}
    \leq M\sum_{v\in V}rk\acute l_{k,l,v}^{t,e}\;\;\;\;\;\;\;\;\;\;\;\;\;\;\;\;\;\;\;\;\;\;\;\;\;\;\;\;\;\;\;\;\;\;\;\;\;\;\;\;\;\;\;\;\;\;
    \nonumber\\
   \forall l \in (\acute{L_{e}{}\cup L_{e}{})},\forall k \in \acute{K_{e}{}}  ,e\in E,t \in T,p \in P,lcv \in LCV
   \nonumber\\
   OR\;\;\forall l \in \acute{L_{e}{}},\forall k \in K_{e}{},e\in E,t \in T,p \in P,lcv \in LCV&& 
  \end{flalign}

  \begin{flalign} \label{cons33}
    \sum_{\acute{t}\leq t\in T}\sum_{lv\in LV}\sum_{rv\in RV}\sum_{v\in V}rkl_{k,l,v}^{lv,rv,\acute{t},e}
    \leq 1\;\;\;\;\;\;\;\;\;\;\;\;\;\;\;\;\;\;\;\;\;\;\;\;\;\;\;\;
    \nonumber\\
    \;\;\;\;\;\;\;\;\;\; \forall l \in (\acute{L_{e}{}\cup L_{e}{})},\forall k \in \acute{K_{e}{}} ,e\in E,t \in T
     \nonumber\\
     OR\;\;\forall l \in \acute{L_{e}{}},\forall k \in K_{e}{} ,e\in E,t \in T&& 
  \end{flalign}

  \begin{flalign} \label{cons34}
    rkl_{k,l,v}^{lv,rv,t,e}=0
     \;\;\;\;\;\;\;\;\;\;\;\;\;\;\;\;\;\;\;\;\;\;\;\;\;\;\;\;\;\;\;\;\;\;\;\;\;\;\;\;\;\;\;\;\;\;\;\;\;\;\;\;\;\;\;\;\;\;\;\;
    \nonumber\\
     \forall l \in (\acute{L_{e}{}\cup L_{e}{})},\forall k \in \acute{K_{e}{}}   ,e\in E,t \in T,lv \in LV,rv \in RV, \nonumber\\v \in V,
    v\neq pipeline\;\;route\;\; transportation\;\;mode
    \nonumber\\
    OR\;\; \forall l \in \acute{L_{e}{}},\forall k \in K_{e}{}  ,e\in E,t \in T,lv \in LV,rv \in RV, \nonumber\\v \in V,
    v\neq pipeline\;\;route\;\; transportation\;\;mode&&
  \end{flalign}

  \begin{flalign} \label{cons35}
    rk\acute{l}_{k,l,v}^{t,e}=0
     \;\;\;\;\;\;\;\;\;\;\;\;\;\;\;\;\;\;\;\;\;\;\;\;\;\;\;\;\;\;\;\;\;\;\;\;\;\;\;\;\;\;\;\;\;\;\;\;\;\;\;\;\;\;\;\;\;\;\;\;\;\;\;\;\;\;\;
    \nonumber\\
     \;\;\;\;\;\;\;\;\;\;\;\;\;\;\;\;\;\;\;\;\;\;\;\;\;\;\;
     \forall l \in (\acute{L_{e}{}\cup L_{e}{})},\forall k \in \acute{K_{e}{}}   ,e\in E,t \in T, \nonumber\\
    v= pipeline\;\;route\;\; transportation\;\;mode
    \nonumber\\
    OR\;\; \forall l \in \acute{L_{e}{}},\forall k \in K_{e}{}  ,e\in E,t \in T, \nonumber\\
    v= pipeline\;\;route\;\; transportation\;\;mode
    &&
  \end{flalign}

   \begin{flalign} \label{cons36}
    \sum_{p\in P}qkl_{p,k,l,v}^{t,e} \leq\sum_{\acute{t}\leq t\in T}\sum_{lv\in LV}
    \sum_{rv\in RV}clv_{lv}^{}rkl_{k,l,v}^{lv,rv,\acute{t},e}+Rkl_{k,l}^{}\nonumber\\(icapkl_{k,l}^{}+\sum_{\acute{t}\leq t\in T}\sum_{ev\in EV}capkl_{k,l}^{ev}ykl_{k,l}^{ev,\acute{t},e})\;\;\;\;\;\;\;\;\;\;\;\;\;\;\;\;\;\;\;\;\;\;\;\;
    \nonumber\\\forall l \in (\acute{L_{e}{}\cup L_{e}{})},k \in (\acute{K_{e}{}\cup K_{e}{})} ,e\in E,
    t \in T,p \in P,v \in V\;\;\;&& 
  \end{flalign}

   \begin{flalign} \label{cons37}
    qlpl_{p,lp,l,v}^{t,e} \leq\sum_{lcv\in LCV}trc_{v}^{lcv}nlpl_{p,lp,l,v}^{lcv,t,e}
    TPP+Rlpl_{lp,l}^{}\;\;\;\;\;\;\;\;\;\;\;\;\nonumber\\
    (icaplpl_{lp,l}^{}+\sum_{\acute{t}\leq t\in T}\sum_{ev\in EV}caplpl_{lp,l}^{ev}ylpl_{lp,l}^{ev,\acute{t},e})+\sum_{\acute{t}\leq t\in T}\sum_{lv\in LV}
    \nonumber\\
    \sum_{rv\in RV}clv_{lv}^{}rlpl_{lp,l,v}^{lv,rv,\acute{t},e}\;\;\;\forall l \in L_{e}\;\; and\;\; lp \in L_{e},e\in E,t \in T
     && 
  \end{flalign}

  \begin{flalign} \label{cons38}
    \sum_{\acute{t}\leq t\in T}\sum_{lv\in LV}\sum_{rv\in RV}\sum_{v\in V}rlpl_{lp,l,v}^{lv,rv,\acute{t},e}+\sum_{v\in V}rlp\acute l_{lp,l,v}^{t,e}
    \leq 1\;\;\;\;\;\;\;\;\;\;\;\;
    \nonumber\\
     \forall l \in L_{e}^{}, lp \in L_{e}{} ,e\in E,t \in T\;\;\;\;\;
     && 
  \end{flalign}

  \begin{flalign} \label{cons39}
    \sum_{\acute{t}\leq t\in T}\sum_{lv\in LV}\sum_{rv\in RV}\sum_{v\in V}rlpl_{lp,l,v}^{lv,rv,\acute{t},e}+Rlpl_{lp,l}^{}
    \leq 1\;\;\;\;\;\;\;\;\;\;\;\;
    \nonumber\\
    \;\;\;\;\;\;\;\;\;\;\;\;\;\;\;\;\;\;\;\;\;\;\;\;\;\;\;\;\;\;\;\;\;\;\;\;\;\; \forall l \in L_{e}{},lp \in L_{e}{}  ,e\in E,t \in T&& 
  \end{flalign}

  \begin{flalign} \label{cons40}
    nlpl_{p,lp,l,v}^{lcv,t,e}
    \leq M\sum_{v\in V}rlp\acute l_{lp,l,v}^{t,e}\;\;\;\;\;\;\;\;\;\;\;\;\;\;\;\;\;\;\;\;\;\;\;\;\;\;\;\;\;\;\;\;\;\;\;\;\;\;\;\;\;\;\;\;\;\;
    \nonumber\\
    \;\;\;\;\;\;\;\;\;\;\;\;\;\;\forall l \in L_{e}{},lp \in L_{e}{}  ,e\in E,t \in T,p \in P,lcv \in LCV&& 
  \end{flalign}

  \begin{flalign} \label{cons41}
    \sum_{\acute{t}\leq t\in T}\sum_{ev\in EV}ylpl_{lp,l}^{ev,\acute{t},e}\leq Rlpl_{lp,l}^{}
     \;\;\;\;\;\;\;\;\;\;\;\;\;\;\;\;\;\;\;\;\;\;\;\;\;\;\;\;\;\;\;\;\;\;
    \nonumber\\
    \;\;\;\;\;\;\;\;\;\;\;\;\;\;\;\;\;\;\;\;\;\;\;\;\;\;\;\;\;\;\;\;\;\;\;\;\;\; \forall l \in L_{e}{},lp \in L_{e}{}  ,e\in E,t \in T&& 
  \end{flalign}

   \begin{flalign} \label{cons42}
    qlpl_{p,lp,l,v}^{t,e} \leq\sum_{lcv\in LCV}trc_{v}^{lcv}nlpl_{p,lp,l,v}^{lcv,t,e}
    TPP+\sum_{\acute{t}\leq t\in T}\sum_{lv\in LV}
    \nonumber\\
    \sum_{rv\in RV}clv_{lv}^{}rlpl_{lp,l,v}^{lv,rv,\acute{t},e}\;\;\;\;\;\;\;\;\;\;\;\;\;\;\;\;\;\;\;\;\;\;\;\;\;\;\;\;\;\;\;\;\;\;\;\;\;\;\;\;\;\;\;\;\;\;\;\;\;\;\;\;\;\;\nonumber\\\forall  l\;and\;lp \in (\acute{L_{e}{}\cup L_{e}{})},e\in E,t \in T,p \in P
     \;\;\;&& 
  \end{flalign}

  \begin{flalign} \label{cons43}
    \sum_{\acute{t}\leq t\in T}\sum_{lv\in LV}\sum_{rv\in RV}\sum_{v\in V}rlpl_{lp,l,v}^{lv,rv,\acute{t},e}+\sum_{v\in V}rlp\acute l_{lp,l,v}^{t,e}
    \leq 1\;\;\;\;\;\;\;\;\;\;\;\;
    \nonumber\\
     \forall l\;and\;lp \in (\acute{L_{e}{}\cup L_{e}{})},e\in E,t \in T\;\;\;\;\;\;\;\;
     && 
  \end{flalign}

  \begin{flalign} \label{cons44}
    \sum_{\acute{t}\leq t\in T}\sum_{lv\in LV}\sum_{rv\in RV}\sum_{v\in V}rlpl_{lp,l,v}^{lv,rv,\acute{t},e}
    \leq 1\;\;\;\;\;\;\;\;\;\;\;\;\;\;\;\;\;\;\;\;\;\;\;\;\;\;
    \nonumber\\
    \;\;\;\;\;\;\;\;\;\;\;\;\;\;\;\;\;\;\;\;\;\;\;\;\;\;\;\;\;\;\; \forall  l\;and\;lp \in (\acute{L_{e}{}\cup L_{e}{})},e\in E,t \in T&& 
  \end{flalign}

  \begin{flalign} \label{cons45}
    nlpl_{p,lp,l,v}^{lcv,t,e}
    \leq M\sum_{v\in V}rlp\acute l_{lp,l,v}^{t,e}\;\;\;\;\;\;\;\;\;\;\;\;\;\;\;\;\;\;\;\;\;\;\;\;\;\;\;\;\;\;\;\;\;\;\;\;\;\;\;\;\;\;\;\;\;\;
    \nonumber\\
    \forall l\;and\;lp \in (\acute{L_{e}{}\cup L_{e}{})},e\in E,t \in T,p \in P,lcv \in LCV\;\;\;\;&& 
  \end{flalign}

  \begin{flalign} \label{cons46}
    rlpl_{lp,l,v}^{lv,rv,t,e}=0
     \;\;\;\;\;\;\;\;\;\;\;\;\;\;\;\;\;\;\;\;\;\;\;\;\;\;\;\;\;\;\;\;\;\;\;\;\;\;\;\;\;\;\;\;\;\;\;\;\;\;\;\;\;\;\;\;\;\;\;\;\;\;
    \nonumber\\
    \;\;\;\;\;\;\;\;\;\;\;\;\; \forall l\;and\;lp \in (\acute{L_{e}{}\cup L_{e}{})} ,e\in E,t \in T,p \in P,v \in V \nonumber\\
    v\neq pipeline\;\;route\;\; transportation\;\;mode&&
  \end{flalign}

  \begin{flalign} \label{cons47}
    rlp\acute{l}_{lp,l,v}^{t,e}=0
     \;\;\;\;\;\;\;\;\;\;\;\;\;\;\;\;\;\;\;\;\;\;\;\;\;\;\;\;\;\;\;\;\;\;\;\;\;\;\;\;\;\;\;\;\;\;\;\;\;\;\;\;\;\;\;\;
    \nonumber\\
    \;\;\;\;\;\;\;\;\;\;\;\;\; \forall l\;and\;lp \in (\acute{L_{e}{}\cup L_{e}{})}  ,e\in E,t \in T,p \in P \nonumber\\
    v= pipeline\;\;route\;\; transportation\;\;mode&&
  \end{flalign}

   \begin{flalign} \label{cons48}
    \;\;\;qlpl_{p,lp,l,v}^{t,e} \leq Rlpl_{lp,l}^{}(icaplpl_{lp,l}^{}+\sum_{\acute{t}\leq t\in T}\sum_{ev\in EV}caplpl_{lp,l}^{ev}\nonumber\\
    ylpl_{lp,l}^{ev,\acute{t},e})+\sum_{\acute{t}\leq t\in T}\sum_{lv\in LV}\sum_{rv\in RV}clv_{lv}^{}rlpl_{lp,l,v}^{lv,rv,\acute{t},e}
    \;\;\;\;\;\;\;\;\;\;\;\;\nonumber\\
    \forall l\;and\;lp \in (\acute{L_{e}{}\cup L_{e}{})},e\in E,t \in T
     \nonumber\\
    v= pipeline\;\;route\;\; transportation&& 
  \end{flalign}

   \begin{flalign} \label{cons49}
    \;\;\;qlm_{p,l,m,v}^{t,e} \leq\sum_{lcv\in LCV}trc_{v}^{lcv}nlm_{p,l,m,v}^{lcv,t,e}
    TPP\;\;\;\;\;\;\;\;\;\;\;\;\;\;\;\;\;\;\;\;\;
    \nonumber\\
    \;\;\;\forall l\;and\;lp \in (\acute{L_{e}{}\cup L_{e}{})},e\in E,t \in T,v \in V,m \in M,p \in P
     && 
  \end{flalign}

   \begin{flalign} \label{cons50}
    \;\;\;\sum_{k\in (K_{e}\cup\acute{K_{e}})}\sum_{l\in (L_{e}\cup\acute{L_{e}})}\sum_{p\in P}nkl_{p,k,l,v}^{lcv,t,e}+\sum_{lp\in (L_{e}\cup\acute{L_{e}})}\sum_{l\in (L_{e}\cup\acute{L_{e}})}\;\;\;\;\;\;
    \nonumber\\
    \sum_{p\in P}nlpl_{p,lp,l,v}^{lcv,t,e}+\sum_{l\in (L_{e}\cup\acute{L_{e}})}\sum_{p\in P}\sum_{m\in M}nlm_{p,l,m,v}^{lcv,t,e}
    \leq nmax_{}^{lcv}\nonumber\\
    \;\;\;\forall e\in E,t \in T,v \in V,lcv \in LCV\;\;\;\;\;\;
     && 
  \end{flalign}

   \begin{flalign} \label{cons51}
    qkl_{p,k,l,v}^{t,e} \leq M\sum_{\acute{t}\leq t\in T}\sum_{el\in EL}
    xl_{l}^{el,\acute{t},e}\;\;\;\;\;\;\;\;\;\;\;\;\;\;\;\;\;\;\;\;\;\;\;\;\;\;\;\nonumber\\\;\;\;\;\;\;\;\;\;\;\;\;\;\;\;\;\;\;\;\;\;\forall l \in \acute{L_{e}{}},k \in (\acute{K_{e}{}\cup K_{e}{})} ,e\in E,
    t \in T\;\;\;&& 
  \end{flalign}

   \begin{flalign} \label{cons52}
    \;\;qkl_{p,k,l,v}^{t,e} \leq M\sum_{\acute{t}\leq t\in K}\sum_{ek\in EK}
    xk_{k}^{ek,\acute{t},e}\;\;\;\;\;\;\;\;\;\;\;\;\;\;\;\;\;\;\;\;\;\;\;\;\;\nonumber\\\;\;\;\;\;\;\;\;\;\;\;\;\;\;\;\;\;\;\;\;\;\forall k \in \acute{K_{e}{}},l \in (\acute{L_{e}{}\cup L_{e}{})} ,e\in E,
    t \in T\;\;\;&& 
  \end{flalign}

   \begin{flalign} \label{cons53}
    qlm_{p,l,m,v}^{t,e} \leq M\sum_{\acute{t}\leq t\in T}\sum_{el\in EL}
    xl_{l}^{el,\acute{t},e}\;\;\;\;\;\;\;\;\;\;\;\;\;\;\;\;\;\;\;\;\;\;\;\;\;\;\;\nonumber\\\;\;\;\;\;\;\;\;\;\;\;\;\;\;\;\;\;\;\;\;\;\forall l \in \acute{L_{e}{}},m\in M,e\in E,
    t \in T\;\;\;&& 
  \end{flalign}       
       
   \begin{flalign} \label{cons54}
    qlpl_{p,lp,l,v}^{t,e} \leq M\sum_{\acute{t}\leq t\in T}\sum_{el\in EL}
    xl_{l}^{el,\acute{t},e}\;\;\;\;\;\;\;\;\;\;\;\;\;\;\;\;\;\;\;\;\;\;\;\;\;\;\;\nonumber\\\;\;\;\;\;\;\;\;\;\;\;\;\;\;\;\;\;\;\;\;\;\forall l \in \acute{L_{e}{}},lp \in \acute{L_{e}{}},e\in E,
    t \in T\;\;\;&& 
  \end{flalign}       
       
   \begin{flalign} \label{cons55}
    qlpl_{p,lp,l,v}^{t,e} \leq M\sum_{\acute{t}\leq t\in T}\sum_{el\in EL}
    xl_{l}^{el,\acute{t},e}\;\;\;\;\;\;\;\;\;\;\;\;\;\;\;\;\;\;\;\;\;\;\;\;\;\;\;\nonumber\\\;\;\;\;\;\;\;\;\;\;\;\;\;\;\;\;\;\;\;\;\;\forall l \in \acute{L_{e}{}},lp \in L_{e}{},e\in E,
    t \in T\;\;\;&& 
  \end{flalign}       
       
  \begin{flalign} \label{cons56}
    \sum_{t\in T}\sum_{ez\in EZ}\sum_{p\in P}
    Nct_{l}^{ez}n_{p,l}^{ez,t,e}\leq\sum_{t\in T}\sum_{el\in EL}
    Ncl_{l}^{el}xl_{l}^{el,t,e} \nonumber\\\forall l \in \acute{L_{e}{}},e\in E,
    t \in T\;\;\;&& 
  \end{flalign}

  \begin{flalign} \label{cons57}
    \sum_{\acute{t}\leq t\in T}\sum_{lv\in LV}\sum_{rv\in RV}rkl_{k,l,v}^{lv,rv,\acute{t},e}\leq \sum_{\acute{t}\leq t\in K}\sum_{ek\in EK}
    xk_{k}^{ek,\acute{t},e}\;\;\;\;\;\;\;\;\;\;\;\;
    \nonumber\\
     \forall k \in \acute{K_{e}{}},l \in (\acute{L_{e}{}\cup L_{e}{})} ,e\in E,
    t \in T\;\;\;\;\;\nonumber\\
    v= pipeline\;\;route\;\; transportation\;\;mode&& 
  \end{flalign}

  \begin{flalign} \label{cons58}
    \sum_{\acute{t}\leq t\in T}\sum_{lv\in LV}\sum_{rv\in RV}rkl_{k,l,v}^{lv,rv,\acute{t},e}\leq \sum_{t\in T}\sum_{el\in EL}
    xl_{l}^{el,t,e}\;\;\;\;\;\;\;\;\;\;\;\;
    \nonumber\\
     \forall k \in \acute{K_{e}{}},l \in (\acute{L_{e}{}\cup L_{e}{})} ,e\in E,
    t \in T\;\;\;\;\;\nonumber\\
    v= pipeline\;\;route\;\; transportation\;\;mode&& 
  \end{flalign}

  \begin{flalign} \label{cons59}
    \sum_{\acute{t}\leq t\in T}\sum_{lv\in LV}\sum_{rv\in RV}rlpl_{lp,l,v}^{lv,rv,\acute{t},e}\leq \sum_{t\in T}\sum_{el\in EL}
    xl_{l}^{el,t,e}\;\;\;\;\;\;\;\;\;\;\;\;
    \nonumber\\
     \forall l \in \acute{L_{e}{}},lp \in  L_{e}{} ,e\in E,
    t \in T\;\;\;\;\;\nonumber\\
    v= pipeline\;\;route\;\; transportation\;\;mode&& 
  \end{flalign}

  \begin{flalign} \label{cons60}
    \sum_{\acute{t}\leq t\in T}\sum_{lv\in LV}\sum_{rv\in RV}rlpl_{lp,l,v}^{lv,rv,\acute{t},e}\leq \sum_{t\in T}\sum_{el\in EL}
    xl_{l}^{el,t,e}\;\;\;\;\;\;\;\;\;\;\;\;
    \nonumber\\
     \forall l \in \acute{L_{e}{}}, lp \in \acute{L_{e}{}} ,e\in E,
    t \in T\;\;\;\;\;\nonumber\\
    v= pipeline\;\;route\;\; transportation\;\;mode&& 
  \end{flalign}

   \begin{flalign} \label{cons61}
    \sum_{k\in (K_{e}\cup\acute{K_{e}})}\sum_{l\in (L_{e}\cup\acute{L_{e}})}\sum_{p\in P}\sum_{v\in V}\sum_{e\in E}qkl_{p,k,l,v}^{t,e}(1\setminus\mu_{p}^{})\leq w_{i}^{}\;\;\;\;\;\;
    \nonumber\\
    \;\;\;\forall i\in I,t \in T\;\;\;\;\;\;
     && 
  \end{flalign}

   \begin{flalign} \label{cons62}
    qlpl_{p,lp,l,v}^{t,e} \leq M
    s_{lp,l}^{t,e}\;\;\;\;\;\;\;\;\;\;\;\;\;\;\;\;\;\;\;\;\;\;\;\;\;\;\;\;\;\;\;\;\;\;\;\;\;\;\;\;\;\;\;\;\;\;\;\;\;\;\;\;\;\;\;\;\;\nonumber\\\;\;\;\;\;\;\;\;\;\;\;\;\;\;\;\;\;\;\;\;\;\forall l \in L_{e}{},lp \in L_{e}{},e\in E,p\in P,v\in V,
    t \in T\;\;\;&& 
  \end{flalign}

   \begin{flalign} \label{cons63}
    qlpl_{p,l,lp,v}^{t,e} \leq M
    s_{l,lp}^{t,e}\;\;\;\;\;\;\;\;\;\;\;\;\;\;\;\;\;\;\;\;\;\;\;\;\;\;\;\;\;\;\;\;\;\;\;\;\;\;\;\;\;\;\;\;\;\;\;\;\;\;\;\;\;\;\;\;\;\nonumber\\\;\;\;\;\;\;\;\;\;\;\;\;\;\;\;\;\;\;\;\;\;\forall l \in L_{e}{},lp \in L_{e}{},e\in E,p\in P,v\in V,
    t \in T\;\;\;&& 
  \end{flalign}

   \begin{flalign} \label{cons64}
    s_{l,lp}^{t,e}+s_{lp,l}^{t,e}\leq 1\;\;\;\;\;\;\;\;\;\;\;\;\;\;\;\;\;\;\;\;\;\;\;\;\;\;\;\;\;\;\;\;\;\;\;\;\;\;\;\;\;\;\;\;\;\;\;\;\;\;\;\;\;\;\;\;\;\nonumber\\\;\;\;\;\;\;\;\;\;\;\;\;\;\;\;\;\;\;\;\;\;\forall l \in L_{e}{},lp \in L_{e}{},t \in T,e\in E\;\;\;&& 
  \end{flalign}

   \begin{flalign} \label{cons65}
    qlpl_{p,l,lp,v}^{t,e} \leq M
    s_{l,lp}^{t,e}\;\;\;\;\;\;\;\;\;\;\;\;\;\;\;\;\;\;\;\;\;\;\;\;\;\;\;\;\;\;\;\;\;\;\;\;\;\;\;\;\;\;\;\;\;\;\;\;\;\;\;\;\;\;\;\;\;\nonumber\\\;\;\;\;\;\;\;\;\;\;\;\;\;\;\;\;\;\;\;\;\;\forall l \in L_{e}{},lp \in L_{e}{},e\in E,p\in P,v\in V,
    t \in T\;\;\;&& 
  \end{flalign}

   \begin{flalign} \label{cons66}
    qlpl_{p,lp,l,v}^{t,e} \leq M
    s_{lp,lk}^{t,e}\;\;\;\;\;\;\;\;\;\;\;\;\;\;\;\;\;\;\;\;\;\;\;\;\;\;\;\;\;\;\;\;\;\;\;\;\;\;\;\;\;\;\;\;\;\;\;\;\;\;\;\;\;\;\;\;\;\nonumber\\\;\;\;\;\;\;\;\;\;\;\;\;\;\;\;\;\;\;\;\;\;\forall lk \in L_{e}{},lp \in L_{e}{},e\in E,p\in P,v\in V,
    t \in T\;\;\;&& 
  \end{flalign}

   \begin{flalign} \label{cons67}
    qlpl_{p,lk,l,v}^{t,e} \leq M
    s_{lk,l}^{t,e}\;\;\;\;\;\;\;\;\;\;\;\;\;\;\;\;\;\;\;\;\;\;\;\;\;\;\;\;\;\;\;\;\;\;\;\;\;\;\;\;\;\;\;\;\;\;\;\;\;\;\;\;\;\;\;\;\;\nonumber\\\;\;\;\;\;\;\;\;\;\;\;\;\;\;\;\;\;\;\;\;\;\forall lk \in L_{e}{},l \in L_{e}{},e\in E,p\in P,v\in V,
    t \in T\;\;\;&& 
  \end{flalign}

   \begin{flalign} \label{cons68}
    s_{l,lp}^{t,e}+s_{lp,lk}^{t,e}+s_{lk,l}^{t,e}\leq 1\;\;\;\;\;\;\;\;\;\;\;\;\;\;\;\;\;\;\;\;\;\;\;\;\;\;\;\;\;\;\;\;\;\;\;\;\;\;\;\;\;\;\;\;\;\;\;\;\;\;\;\;\;\;\;\;\;\nonumber\\\;\;\;\;\;\;\;\;\;\;\;\;\;\;\;\;\;\forall l \in L_{e}{},lp \in L_{e}{},lk \in L_{e}{},t \in T,e\in E\;\;\;\;\;\;&& 
  \end{flalign}

   \begin{flalign} \label{cons69}
   \sum_{l\in (L_{e}\cup\acute{L_{e}})}\sum_{v\in V}\sum_{e\in E}qlm_{p,l,m,v}^{t,e}\geq d_{p,m}^{t}\;\;\;\;\;\;
    \nonumber\\
    \;\;\;\;\;\;\;\;\;\forall p\in P,t \in T,m\in M
     && 
  \end{flalign}

   \begin{flalign} \label{cons70}
   \sum_{l\in (L_{e}\cup\acute{L_{e}})}\sum_{v\in V}\sum_{m\in M}qlm_{p,l,m,v}^{t,e}\geq D_{p}^{t,e}\;\;\;\;\;\;
    \nonumber\\
    \;\;\;\;\;\;\;\;\;\forall p\in P,t \in T,e\in E
     && 
  \end{flalign}

   \begin{flalign} \label{cons71}
   \sum_{lp\in (L_{e}\cup\acute{L_{e}})}\sum_{v\in V}qlpl_{p,lp,l,v}^{t,e}-\sum_{l\in \acute{ L_{e}}}\sum_{v\in V}qlpl_{p,l,lp,v}^{t,e}+\sum_{k\in (K_{e}\cup\acute{K_{e}})}
   \nonumber\\\sum_{v\in V}qkl_{p,k,l,v}^{t,e}+ vl_{p,l}^{t-1,e}+ i_{p,l}^{t,e}=\sum_{v\in V}\sum_{m\in M}qlm_{p,l,m,v}^{t,e}\;\;\;\;\;\;\;\;\;\nonumber\\
   +Ep_{p,l}^{t,e}+vl_{p,l}^{t,e}
   \;\;\;\;\;\;\;\;\;\forall p\in P,t \in T,e\in E,l\in \acute{ L_{e}}
    \;\;\;\;\;\;\;\;\;\;\;\;\;\;\; && 
  \end{flalign}

   \begin{flalign} \label{cons72}
   \sum_{lp\in L_{e}}\sum_{v\in V}qlpl_{p,lp,l,v}^{t,e}-\sum_{lp\in (L_{e}\cup\acute{L_{e}})}\sum_{v\in V}qlpl_{p,l,lp,v}^{t,e}+\sum_{k\in (K_{e}\cup\acute{K_{e}})}
   \nonumber\\\sum_{v\in V}qkl_{p,k,l,v}^{t,e}+ vl_{p,l}^{t-1,e}+ i_{p,l}^{t,e}=\sum_{v\in V}\sum_{m\in M}qlm_{p,l,m,v}^{t,e}\;\;\;\;\;\;\;\;\;\nonumber\\
   +Ep_{p,l}^{t,e}+vl_{p,l}^{t,e}
   \;\;\;\;\;\;\;\;\;\forall p\in P,t \in T,e\in E,l\in  L_{e}
    \;\;\;\;\;\;\;\;\;\;\;\;\;\;\; && 
  \end{flalign}

   \begin{flalign} \label{cons73}
   vk_{k}^{t,e}\leq ick_{k}^{}(1-\sum_{\acute{t}\leq t\in T}\psi k_{k}^{t,e})+\sum_{\acute{t}\leq t\in T}
   \sum_{uk\in UK}capk_{k}^{uk}
   \nonumber\\+ \tau k_{k}^{uk,\acute{t},e} \;\;\;\;\;\;\;\;\;\forall t \in T,e\in E,k\in  K_{e}
    \;\;\;\;\;\;\;\;\;\;\;\;\;\;\; && 
  \end{flalign}

  \begin{flalign} \label{cons74}
   \;\;\;\sum_{e\in E}vl_{p,l}^{t,e}\leq \sum_{e\in E}ic_{l}^{p}\xi_{l}^{t,e}+\sum_{\acute{t}\leq t\in T}
   \sum_{ul\in UL}\sum_{e\in E}capl_{l}^{ul} 
   \nonumber\\+ \tau l_{p,l}^{ul,\acute{t},e} \;\;\;\;\;\;\;\;\;\forall t \in T,p\in P,l\in  L_{e}
    \;\;\;\;\;\;\;\;\;\;\;\;\;\;\; && 
  \end{flalign}

  \begin{flalign} \label{cons75}
  \;\;\;vk_{k}^{t,e}\leq \sum_{\acute{t}\leq t\in T}
   \sum_{ek\in EK}Nck_{k}^{ek}xk_{k}^{ek,\acute{t},e}
   \nonumber\\  \;\;\;\;\;\;\;\;\;\;\;\;\;\;\;\;\;\;\;\;\;\;\;\;\;\;\;\forall t \in T,k\in  \acute{K_{e}}
     && 
  \end{flalign}

  \begin{flalign} \label{cons76}
  \;\;\;vl_{p,l}^{t,e}\leq \sum_{\acute{t}\leq t\in T}
   \sum_{ez\in EZ}Nct_{l}^{ez}n_{p,l}^{ez,\acute{t},e}
   \nonumber\\  \;\;\;\;\;\;\;\;\;\;\;\;\forall t \in T,p \in P,e \in E,l\in  \acute{L_{e}}
     && 
  \end{flalign}

   \begin{flalign} \label{cons77}
   vk_{k}^{t,e}\geq lk_{k}^{}(ick_{k}^{}(1-\sum_{\acute{t}\leq t\in T}\psi k_{k}^{t,e})+\sum_{\acute{t}\leq t\in T}
   \sum_{uk\in UK}capk_{k}^{uk}
   \nonumber\\+ \tau k_{k}^{uk,\acute{t},e}) \;\;\;\;\;\;\;\;\;\forall t \in T,e\in E,k\in  K_{e}
    \;\;\;\;\;\;\;\;\;\;\;\;\;\;\; && 
  \end{flalign}

  \begin{flalign} \label{cons78}
   \;\;\;\sum_{e\in E}vl_{p,l}^{t,e}\geq l_{l}^{}(\sum_{e\in E}ic_{l}^{p}\xi_{l}^{t,e}+\sum_{\acute{t}\leq t\in T}
   \sum_{ul\in UL}\sum_{e\in E}capl_{l}^{ul} 
   \nonumber\\+ \tau l_{p,l}^{ul,\acute{t},e}) \;\;\;\;\;\;\;\;\;\forall t \in T,p\in P,l\in  L_{e}
    \;\;\;\;\;\;\;\;\;\;\;\;\;\;\; && 
  \end{flalign}

  \begin{flalign} \label{cons79}
  \;\;\;vk_{k}^{t,e}\geq lk_{k}^{} (\sum_{\acute{t}\leq t\in T}
   \sum_{ek\in EK}Nck_{k}^{ek}xk_{k}^{ek,\acute{t},e})
   \nonumber\\  \;\;\;\;\;\;\;\;\;\;\;\;\;\;\;\;\;\;\;\;\;\;\;\;\;\;\;\forall t \in T,e \in E,k\in  \acute{K_{e}}
     && 
  \end{flalign}

  \begin{flalign} \label{cons80}
  \;\;\;vl_{p,l}^{t,e}\geq l_{l}^{}( \sum_{\acute{t}\leq t\in T}
   \sum_{ez\in EZ}Nct_{l}^{ez}n_{p,l}^{ez,\acute{t},e})
   \nonumber\\  \;\;\;\;\;\;\;\;\;\;\;\;\forall t \in T,p \in P,e \in E,l\in  \acute{L_{e}}
     && 
  \end{flalign}

According to inequalities \ref{cons6} and \ref{cons7} , at most one facility (refinery or distribution center) at a given capacity level is selected by a stakeholder in each candidate site. The number of storage tanks for each refined product (p) in each distribution center is determined only once in the model, as shown by relations \ref{cons8}-\ref{cons9}. Constraints \ref{cons10}  and \ref{cons11}  show that only one capacity expansion with a known capacity level and only by one stakeholder takes place in each refinery or distribution center, respectively. Similarly, considering constraints \ref{cons12} and \ref{cons13}, for each pipeline from refineries to distribution centers and between distribution centers, respectively, only one capacity expansion with a known capacity level and only by one stakeholder takes place. Constraint \ref{cons14} shows that a refinery closes down only once during modeling periods. It must be pointed out that, once the refinery was closed down it may not be opened again. As stated in the problem definition section, distribution centers are shared sources and are simultaneously used by several stakeholders. $\xi_{l}^{t,e} $indicates the share of stakeholder e in the distribution center l during time period t and is less than one for each stakeholder according to constraint \ref{cons15}. Constraint \ref{cons16} shows that the sum of stakeholders’ shares in a distribution center must be less than one. Constraint \ref{cons17} shows if a distribution center is decided to be closed down within time periods, it should not be expanded. Input streams to existing refineries and distribution centers are limited to their current capacity and the capacity expansions is performed on them based on equations \ref{cons18} and \ref{cons19}. According to inequalities \ref{cons20} and \ref{cons21}, the stream between refineries and distribution centers must fall within a certain range (between the defined minimum and maximum amounts of stream). Therefore, a refinery or distribution center is constructed if the stream from the refineries to the distribution centers is greater than minimum defined stream to each refinery and distribution center. In addition, if and only if a distribution center is constructed, the refined products are stored in it, which is shown by Equation \ref{cons22}. The transportation of refined products between refineries and distribution centers takes place by one of the transportation modes of pipeline, road or rail. Constraints \ref{cons23}-\ref{cons29} show the streams of refined products and selection of transportation modes between existing refineries and distribution centers. Constraints \ref{cons30}-\ref{cons35} show the streams of refined products and selection of transportation modes between new refineries and distribution centers. Constraint \ref{cons36} shows that the sum of the product streams passing through a pipeline must be smaller than its capacity. Relations \ref{cons37}-\ref{cons47} are similar to constraints \ref{cons23}-\ref{cons35} and are used to determine the streams of refined products and select the transportation mode between two distribution centers. Similar to constraint \ref{cons36}, constraint \ref{cons48} shows that the sum of the product streams passing through a pipeline must be smaller than its capacity. The refined products are transported from the distribution centers to customers via road and rail modes, but not pipelines. The selection of transportation mode between distribution centers and customers is determined by relation \ref{cons49}.

Constraint \ref{cons50} shows transportation mode capacity. The maximum number of available facilities from each transportation mode is shown by equation \ref{cons50}.  Constraints \ref{cons51}-\ref{cons55} show that the stream between a refinery and distribution center exits if and only if new refineries and distribution centers are constructed. According to inequality \ref{cons56}, the sum of storage tank capacities in each distribution center must be lower than its capacity. Constraints \ref{cons57}-\ref{cons60} show that a new pipeline to or from a new facility is constructed if and only if those new facilities are constructed. Constraint \ref{cons61} shows the limitation of oil stream input to refineries to the crude oil extracted from oil fields. Constraints \ref{cons62}-\ref{cons68} have been defined to prevent loop formation. Considering constraint \ref{cons69}, the demand of each area customer for each product must be supplied by all the stakeholders. The amount of refined products supplied to the final customers by each stakeholder must be equal to the demand planned for each stakeholder, which is shown by inequality \ref{cons70}. Equations \ref{cons71} and \ref{cons72} show that the amount of input stream to each facility must be equal to the corresponding output. The inventory stored in each facility must be less than the storage capacity of the facility, which is shown by Equations \ref{cons73}-\ref{cons76}. In addition, the inventory stored in each facility must be greater than the defined minimum inventory, which is shown by Equations \ref{cons77}-\ref{cons80}.

\subsubsection{Job creation constraints}

  \begin{flalign} \label{cons81}
  RLab_{en}^{lev,t}\leq ALab_{en}^{lev,t}\;\;\;\;\;\;\;\;\;\;\;\;\;\;\;\;\;\;\;\;\;\;\;\;\;\;\;\;\;\;\;\;\;\;\;\;
   \nonumber\\  \;\;\;\;\;\;\;\;\;\;\;\;\;\;\;\;\;\;\;\;\;\;\;\;\;\;\;\forall t \in T,en \in EN,lev \in LEV
     && 
  \end{flalign}

  \begin{flalign} \label{cons82}
  ALab_{en}^{lev,t}=ALab_{en}^{lev,t-1}- RLab_{en}^{lev,t-1}+NLab_{en}^{lev,t}\;\;\;\;\;\;
   \nonumber\\  \;\;\;\;\;\;\;\;\;\;\;\;\;\;\;\;\;\;\;\;\;\;\;\;\;\;\;\forall t \in T,en \in EN,lev \in LEV
     && 
  \end{flalign}

  \begin{flalign} \label{cons83}
  \sum_{\acute{en}\in EN}\sum_{k\in \acute{K_{e}}}HENK_{en,k}^{\acute{en},lev,t}+\sum_{\acute{en}\in EN}\sum_{k\in \acute{K_{e}}}HEEK_{en,k}^{\acute{en},lev,t}
  \nonumber\\ +\sum_{\acute{en}\in EN}\sum_{l\in \acute{L_{e}}}HENL_{en,l}^{\acute{en},lev,t}+\sum_{\acute{en}\in EN}\sum_{l\in \acute{L_{e}}}HEEL_{en,l}^{\acute{en},lev,t}
  \nonumber\\ =RLab_{en}^{lev,t}\;\;\;\;\;\;\;\;\;\;\;\;\;\;\;\;\;\;\;\;\;\;\;\;\;\;\;\;\;\;\;\;\;\;\;\;\;\;\;\;\;\;\;\;\;\;\;\;\;\;\;\;\;\;\;\;\;\;\;\;\;\;\;\;\;\;
   \nonumber\\  \;\;\;\;\;\;\;\;\;\;\;\;\;\;\;\;\;\;\;\;\;\;\;\;\;\;\;\forall t \in T,en \in EN,lev \in LEV
     && 
  \end{flalign}

  \begin{flalign} \label{cons84}
  \sum_{ek\in EK}Nck_{k}^{ek}xk_{k}^{ek,t,e}WNK_{lev}^{}=\sum_{\acute{en}\in EN}\sum_{en\in EN}Nk_{en}^{k}\;\;\;
  \nonumber\\  HENK_{en,k}^{\acute{en},lev,t}\;\;\;\;\;\;\;\;\forall t \in T,e \in E,lev \in LEV,k\in \acute{K_{e}}
     && 
  \end{flalign}

  \begin{flalign} \label{cons85}
  \sum_{uk\in UK}capk_{k}^{uk}\tau k_{k}^{uk,t,e}WEK_{lev}^{}=\sum_{\acute{en}\in EN}\sum_{en\in EN}Nek_{en}^{k}\;\;\;
  \nonumber\\  HEEK_{en,k}^{\acute{en},lev,t}\;\;\;\;\;\;\;\;\forall t \in T,e \in E,lev \in LEV,k\in K_{e}
     && 
  \end{flalign}

  \begin{flalign} \label{cons86}
  \sum_{el\in EL}Ncl_{l}^{el}xl_{l}^{el,t,e}WNL_{lev}^{}=\sum_{\acute{en}\in EN}\sum_{en\in EN}Nl_{en}^{l}\;\;\;\;\;\;\;\;\;
  \nonumber\\  W_{en}^{\acute{en}}HENl_{en,l}^{\acute{en},lev,t}\;\;\;\;\;\;\;\;\forall t \in T,e \in E,lev \in LEV,l\in \acute{L_{e}}
     && 
  \end{flalign}

  \begin{flalign} \label{cons87}
  \sum_{ul\in UL}\sum_{p\in P}capl_{l}^{ul}\tau l_{l}^{ul,t,e}WEL_{lev}^{}=\sum_{\acute{en}\in EN}\sum_{en\in EN}Nel_{en}^{l}\;\;\;
  \nonumber\\  HEEL_{en,l}^{\acute{en},lev,t}\;\;\;\;\;\;\;\;\forall t \in T,e \in E,lev \in LEV,k\in K_{e}
     && 
  \end{flalign}

  According to Equation \ref{cons81}, the number of required labor force in each time period must be smaller than the number of available people. The number of available labor force in a time period is obtained by subtracting the number of required labor force from the number of available labor force in the previous period and adding the number of new standby labor force, which is shown in Equation \ref{cons82}.  The number of required labor force in each time period, which is the number of labor force employed in the refineries and distribution centers, is obtained from Equation \ref{cons83}. According to Equations \ref{cons84}-\ref{cons87}, the number of required labor force in each refinery and distribution center must be equal to the number of employees. Equation \ref{cons88} is referred to as coverage constraint, based on which the sum of constructions, capacity expansions and existing facilities in each area must be less than the minimum defined for that area.

  \subsubsection{Coverage constraint}

  \begin{flalign} \label{cons88}
 \sum_{t\in T}\sum_{ek\in EK}\sum_{e\in E}\sum_{k\in \acute{K_{e}^{}}}xk_{k}^{ek,t,e}
  Nk_{en}^{k}+\sum_{t\in T}
  \sum_{uk\in UK}
  \sum_{e\in E}
  \sum_{k\in K_{e}^{}}
   \nonumber\\ 
   \tau k_{k}^{uk,t,e}Nek_{en}^{k}+\sum_{t\in T}\sum_{el\in EL}\sum_{e\in E}\sum_{l\in \acute{L_{e}^{}}}xl_{l}^{el,t,e}Nl_{en}^{l}+\sum_{t\in T}
  \nonumber\\ 
  \sum_{ul\in UL}\sum_{e\in E}
  \sum_{l\in L_{e}^{}}\tau l_{l}^{ul,t,e}Nel_{en}^{l}NR_{en}^{}+ND_{en}^{}\leq Maxnum_{en}^{}
  \nonumber\\ \;\;\;\;\forall t \in T,en \in EN
     && 
  \end{flalign}

\section{Fuzzy programming technique for a multi-objective, linear problem}\label{sec:Fuzzyprog}

Bellman and Zadeh proposed fuzzy theory for uncertainty management for the first time \cite{bellman1970decision}. The concept of fuzzy theory with appropriate objective functions was used by Zimmermann for solving multi objective, linear programming.  This method can consider minimization of some objective functions along with maximization of other objective functions to establish a balance between them \cite{zimmermann1978fuzzy} . In other words, fuzzy approach can be used to obtain unique responses for multi objective optimization problems. Consider a mathematical programming problem with two objective functions such as Equation \ref{cons89}:

  \begin{flalign} \label{cons89}
 max\;\;f(x)\;\;\;\;\;\;\;\;\;\;\;\;\;\;\;\;\;\;\;\;\;\;\;\;\;\;\;\;\;\;\;\;\;\;\;\;\;\;\;\;\;\;\;\;\;\;\;\;\;\;\;\;\;\;\;\;\;\;\;\;
   \nonumber\\ 
   min\;\;g(x)\;\;\;\;\;\;\;\;\;\;\;\;\;\;\;\;\;\;\;\;\;\;\;\;\;\;\;\;\;\;\;\;\;\;\;\;\;\;\;\;\;\;\;\;\;\;\;\;\;\;\;\;\;\;\;\;\;\;\;\;\;
  \nonumber\\ 
  \;\;\;\;\;\;subject\;\;to:\;\;\;\;\;\;\;\;\;\;\;\;\;\;\;\;\;\;\;\;\;\;\;\;\;\;\;\;\;\;\;\;\;\;\;\;\;\;\;\;\;\;\;\;\;\;\;\;\;\;\;\;\;\;\;\;\;\;
  \nonumber\\ \;\;\;\;h(x)\leq0\;\;\;\;\;\;\;\;\;\;\;\;\;\;\;\;\;\;\;\;\;\;\;\;\;\;\;\;\;\;\;\;\;\;\;\;\;\;\;\;\;\;\;\;\;\;\;\;\;
     && 
  \end{flalign}

The steps corresponding to fuzzy theory for solving a multi objective problem are as follows:

Step 1: Only one objective function is maximized with all constraints at each time, as follows:

  \begin{flalign} \label{con90}
 max\;\;f(x)\;\;\;\;\;\;\;\;\;\;\;\;\;\;\;\;\;\;\;\;\;\;\;\;\;\;\;\;\;\;\;\;\;\;\;max\;\;g(x)\;\;\;\;\;\;\;
  \nonumber\\ 
  \;\;subject\;\;to:\;\;\;\;\;\;\;\;\;\;\;\;\;\;\;\;\;\;\;\;\;\;\;\;\;\;\;\;\;\;\;subject\;\;to:\;\;\;\;\;\;\;
  \nonumber\\ \;\;\;\;h(x)\leq0\;\;\;\;\;\;\;\;\;\;\;\;\;\;\;\;\;\;\;\;\;\;\;\;\;\;\;\;\;\;\;\;\;h(x)\leq0
     && 
  \end{flalign} 
  
The results of solving models are called $f_{}^{*}$  and $g_{}^{*}$.

Step 2: Only one objective function is minimized with all constraint at each time, as follows:

  \begin{flalign} \label{con91}
 min\;\;f(x)\;\;\;\;\;\;\;\;\;\;\;\;\;\;\;\;\;\;\;\;\;\;\;\;\;\;\;\;\;\;\;\;\;\;\;min\;\;g(x)\;\;\;\;\;\;\;
  \nonumber\\ 
  \;\;subject\;\;to:\;\;\;\;\;\;\;\;\;\;\;\;\;\;\;\;\;\;\;\;\;\;\;\;\;\;\;\;\;\;\;subject\;\;to:\;\;\;\;\;\;\;
  \nonumber\\ \;\;\;\;h(x)\leq0\;\;\;\;\;\;\;\;\;\;\;\;\;\;\;\;\;\;\;\;\;\;\;\;\;\;\;\;\;\;\;\;\;h(x)\leq0
     && 
  \end{flalign}

The results of solving models are called $f_{}^{-}$  and $g_{}^{-}$.

Step 3: Finally, the objective functions are obtained using Equations \ref{cons92} and \ref{cons93} based on f(x) and g(x) functions, which are maximization and minimization objective functions, respectively.

  \[\label{cons92}
    \mu(f)=\left\{
                \begin{array}{ll}
                  0\;\;\;\;\;\;\;\;\;\;\;\;\;\;\;\;\;\;\;\;\;f(x)\leq f_{}^{-}\\
                 \frac{f(x)-f_{}^{-}}{f_{}^{*}-f_{}^{-}}\;\;\;\;\;\;\;\;\;\;f_{}^{-}\leq f(x)\leq f_{}^{*}\\
                  1\;\;\;\;\;\;\;\;\;\;\;\;\;\;\;\;\;\;\;\;\;f(x)\geq f_{}^{*}
                \end{array}
              \right.
  \]

  \[\label{cons93}
    \mu(g)=\left\{
                \begin{array}{ll}
                  0\;\;\;\;\;\;\;\;\;\;\;\;\;\;\;\;\;\;\;\;\;g(x)\leq g_{}^{-}\\
                 \frac{g_{}^{*}-g(x)}{g_{}^{*}-g_{}^{-}}\;\;\;\;\;\;\;\;\;\;g_{}^{-}\leq g(x)\leq g_{}^{*}\\
                  1\;\;\;\;\;\;\;\;\;\;\;\;\;\;\;\;\;\;\;\;\;g(x)\geq g_{}^{*}
                \end{array}
              \right.
  \]

These objective functions must be greater than $\lambda$. Therefore, the multi objective, linear programming model is modeled using fuzzy programming as follows:

  \begin{flalign} \label{cons94}
 max\;\;\lambda\;\;\;\;\;\;\;\;\;\;\;\;\;\;\;\;\;\;\;\;\;\;\;\;\;\;\;\;\;\;\;\;\;\;\;\;\;\;\;\;\;\;\;\;\;\;\;\;\;\;\;\;\;\;\;\;\;\;\;\;\;\;\;\;\;\;\;
 \nonumber\\ 
  \;\;\;\;\;\;subject\;\;to:\;\;\;\;\;\;\;\;\;\;\;\;\;\;\;\;\;\;\;\;\;\;\;\;\;\;\;\;\;\;\;\;\;\;\;\;\;\;\;\;\;\;\;\;\;\;\;\;\;\;\;\;\;\;\;\;\;\;\;\;
  \nonumber\\ \;\;\;\;h(x)\leq0\;\;\;\;\;\;\;\;\;\;\;\;\;\;\;\;\;\;\;\;\;\;\;\;\;\;\;\;\;\;\;\;\;\;\;\;\;\;\;\;\;\;\;\;\;\;\;\;\;
  \nonumber\\ \;\;\;\;f(x)\geq f_{}^{-}+\lambda(f_{}^{*}-f_{}^{-})\;\;\;\;\;\;\;\;\;\;\;\;\;\;\;\;\;\;\;\;\;\;\;\;\;
  \nonumber\\ \;\;\;\;g(x)\leq g_{}^{*}+\lambda(g_{}^{*}-g_{}^{-})\;\;\;\;\;\;\;\;\;\;\;\;\;\;\;\;\;\;\;\;\;\;\;\;\;\;
     && 
  \end{flalign} 
 
 Multi- period, multi objective, multi-stakeholder, multi-level and multi-product DSCPSC model is converted into a single-objective problem via fuzzy approach using the above linear programming model.

\section{Numerical results}\label{sec:Numericalresult}

As stated, the multi objective, multi-stakeholder, multi-level and multi-product DSCPSC model is an MILP model in which the outputs are flow rates, installation of new facilities, expansion of existing facilities, installation of new pipelines, expansion of existing pipelines, exports, imports, determination of each of the stakeholder’s share of distribution centers, inventory and selection of transportation modes for each of the stakeholders. Three stakeholders have been considered in the model.

\begin{figure}[h]
    \centering
    \includegraphics[width=1.03\linewidth,trim={0.2in 0.2in 0.2in 0.2in},clip]{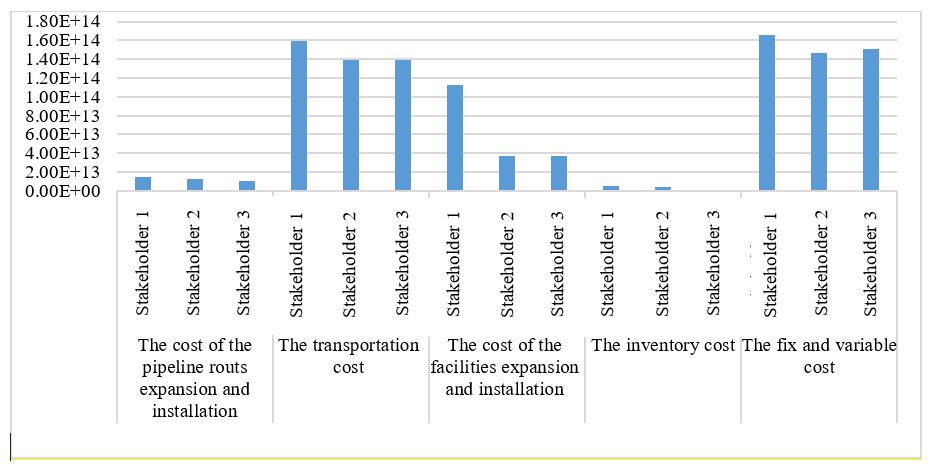}
    \caption{Stakeholders’ costs}
    \label{fig:fig1}
\end{figure}

Since the proposed model is a multi-objective model, fuzzy programming technique is used to solve and convert it to a single-objective model. Therefore, three objective functions are converted into one and three constraints are added to the model constraints. The model is thus solved using GAMS 24.1.2 software. In order to use fuzzy planning technique, the maximum and minimum values of each function are calculated. The maximum and minimum values are then placed in constraint 94.
Therefore, optimal responses are obtained by considering one objective function. Model parameter values are shown in Appendix A..
Model objective function values are shown in Table \ref{table3}.
The corresponding costs for each of the stakeholders including pipeline construction and expansion, facility construction and expansion, transportation of products, inventory and fixed and variable costs in each of the equilibria are shown in Figure \ref{fig:fig1}.

\begin{table}[H]
   
\caption{Optimal values of the objective function for the whole supply chain and separate stakeholders}
\label{table3}
\centering
\begin{tabular}{|>{\centering}c |>{\centering}c |>{\centering}c |}

\hline
 
\hline
  Beneficiary&	Objective function&	Nash game value \tabularnewline \hline

\multirow{3}{*}{Stakeholder 1}  &$P_{e}^{}$&1.08E+15\tabularnewline \cline{2-3}
 &$Pul_{e}^{}$&5.19926E+12	 \tabularnewline \cline{2-3}
  &$S_{e}^{}$&21570\tabularnewline \hline

\multirow{3}{*}{Stakeholder 2}  &$P_{e}^{}$&7.23E+14\tabularnewline \cline{2-3}
 &$Pul_{e}^{}$&3.48061E+12	 \tabularnewline \cline{2-3}
  &$S_{e}^{}$&17989\tabularnewline \hline

\multirow{3}{*}{Stakeholder 3}  &$P_{e}^{}$&8.82E+14\tabularnewline \cline{2-3}
 &$Pul_{e}^{}$&4.24606E+12	 \tabularnewline \cline{2-3}
  &$S_{e}^{}$&19173\tabularnewline \hline

\multicolumn{2}{|c|}{$\lambda$}&0.60105 \tabularnewline \hline

\end{tabular}
\end{table}


As observed, the corresponding pipeline and facility construction and expansion, inventory and fixed and variable costs for the first stakeholder are greater compared with those of the second and third stakeholders. It is also observed that fixed and variable, transportation and facility construction and expansion costs are remarkable in comparison with other costs.

\begin{figure}[h]
    \centering
    \includegraphics[width=0.7\linewidth,trim={0.2in 0.2in 0.2in 0.2in},clip]{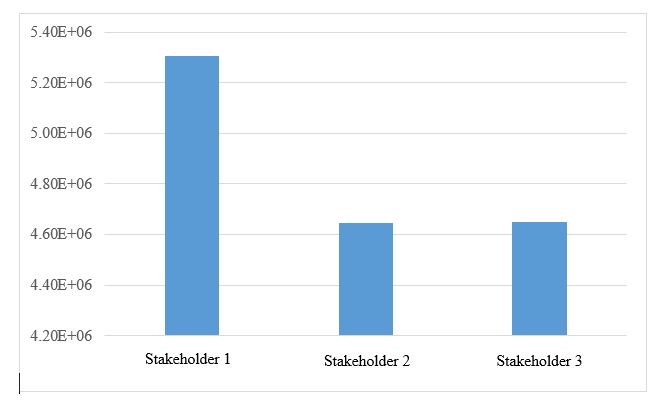}
    \caption{Flow of products transported by stakeholders}
    \label{fig:fig2}
\end{figure}

\begin{table*}[ht]
   
\caption{Sensitivity analysis of changes in the objective function coefficients}
\label{table4}
\centering
\begin{tabular}{|>{\centering}p{1cm} |>{\centering}p{1.1cm} |>{\centering}p{1.5cm} |>{\centering}p{1.5cm} |>{\centering}p{1.5cm} |>{\centering}p{1.1cm} |>{\centering}p{1cm} |>{\centering}p{1.5cm} |>{\centering}p{1.5cm} |>{\centering}p{1.5cm} |}

\hline
 Parameter&	change in parameter\%	&change in objective function value\% (Stakeholder 1)&	change in objective function value\% (Stakeholder 2)&	change in objective function value\% (Stakeholder 3)&	Parameter&	change in parameter\%	&change in objective function value\% (Stakeholder 1)&	change in objective function value\% (Stakeholder 2)	&change in objective function value\% (Stakeholder 3) \tabularnewline \hline

\multirow{6}{1cm}{Pipeline transport 
costs} &30	&0.31&	0.62&	0.58&\multirow{6}{1.4cm}{Costs of other transportation modes except pipeline} &30&	0&	1.26&0    \tabularnewline\cline{2-5}
\cline{7-10}
&20&	0.20&	0.40&	0.39&& 20	&0&	0.84&	0
\tabularnewline\cline{2-5}
\cline{7-10}

&10&	0.10&	0.20&	0.19&& 10&	0&	0.42&	0
\tabularnewline\cline{2-5}
\cline{7-10}

&-10&	-0.10&	-0.20&	-0.19&&-10&	0&	-0.42&	0
\tabularnewline\cline{2-5}
\cline{7-10}

&-20&	-0.20&	-0.40&	-0.39&& -20&	0&	-0.84&	0
\tabularnewline\cline{2-5}\cline{7-10}

&-30&	-0.31	&-0.60&	-0.56&&-30&	0	&-1.26&	0
\tabularnewline\hline

\multirow{6}{1cm}{Facilities’ expansion costs} &30	&0&	0&	0&\multirow{6}{1.4cm}{Inventory costs} &30&0.10&	0.41&	0.13    \tabularnewline\cline{2-5}
\cline{7-10}
&20&	0&	0&0&& 20	&0.07&	0.27&	0.09
\tabularnewline\cline{2-5}
\cline{7-10}

&10&	0&	0&0&& 10&	0.03&	0.13&	0.04
\tabularnewline\cline{2-5}
\cline{7-10}

&-10&	0&	0&0&&-10&	-0.06&	-0.13	&-0.01
\tabularnewline\cline{2-5}
\cline{7-10}

&-20&	0&	0&0&& -20&	-0.12&	-0.28&	-0.04
\tabularnewline\cline{2-5}\cline{7-10}

&-30&	-0.09&	0.15&	0.02&&-30&	-0.21&	-0.47&	0
\tabularnewline\hline

\multirow{6}{1cm}{Variable costs} &30	&2.66&	4.18&	4.02&\multirow{6}{1.4cm}{Routs’ expansion and installation costs} &30&0.11&	0.10&	0.16    \tabularnewline\cline{2-5}
\cline{7-10}
&20&	1.80&	2.69&	2.67&& 20&	0.08&	0.06&	0.11
\tabularnewline\cline{2-5}
\cline{7-10}

&10&	0.90&	1.31&	1.35&& 10&	0.04&	0.03&	0.06
\tabularnewline\cline{2-5}
\cline{7-10}

&-10&	-0.92&	-1.26&	-1.32&&-10&	-0.04&	-0.03&	-0.05
\tabularnewline\cline{2-5}
\cline{7-10}

&-20&	-1.82&	-2.53&	-2.65&& -20&	-0.07&	-0.06&	-0.11
\tabularnewline\cline{2-5}\cline{7-10}

&-30&	-2.64&	-3.95&	-4.02&&-30&	-0.11&	-0.10&	-0.16
\tabularnewline\hline

\multirow{6}{1cm}{Facilities’ installation costs} &30	&1.41&	1.08&	1.09&\multirow{6}{1.4cm}{Fixed cost} &30&0.03&	0.06	&0.06    \tabularnewline\cline{2-5}
\cline{7-10}
&20&	0.94&	0.72&	0.73&& 20&	0.02&	0.04&	0.04
\tabularnewline\cline{2-5}
\cline{7-10}

&10&	0.47&	0.36&	0.37&& 10&	0.01&	0.02&	0.02
\tabularnewline\cline{2-5}
\cline{7-10}

&-10&	-0.47&	-0.35&	-0.36&&-10&	-0.01	&-0.02&	-0.02
\tabularnewline\cline{2-5}
\cline{7-10}

&-20&	-0.94&	-0.72&	-0.73&& -20&	-0.02&	-0.04&	-0.03
\tabularnewline\cline{2-5}\cline{7-10}

&-30&	-1.41&	-1.07&	-1.08&&-30&	-0.11&	-0.10&	-0.16
\tabularnewline\hline
\end{tabular}
\end{table*}


According to Figure \ref{fig:fig2}, the amounts of products transported by the first stakeholder are greater than those by the second and third stakeholders. This verifies the differences between transportation costs in different stakeholders in Figure \ref{fig:fig1}.

Sensitivity analysis is performed for further study. The objective of performing sensitivity analysis is the investigation of analysis results and the analysis of their variations to parameter changes. The variations in results with objective function coefficient changes (cost parameters) have been considered in sensitivity analysis.

Given that parameter changes individually lead to small changes in the objective function in this problem, the parameters are divided into 8 groups and the variations in the objective function with changes in each group of parameters have been studied (Table \ref{table4}).  Each of these examples studies the variations of a group of parameters. The comparison of the values of the objective function in each of these examples with the corresponding optimal values is shown in Table \ref{table4}.

\begin{figure}[h]
    \centering
    \includegraphics[width=1\linewidth,trim={0.2in 0.2in 0.2in 0.2in},clip]{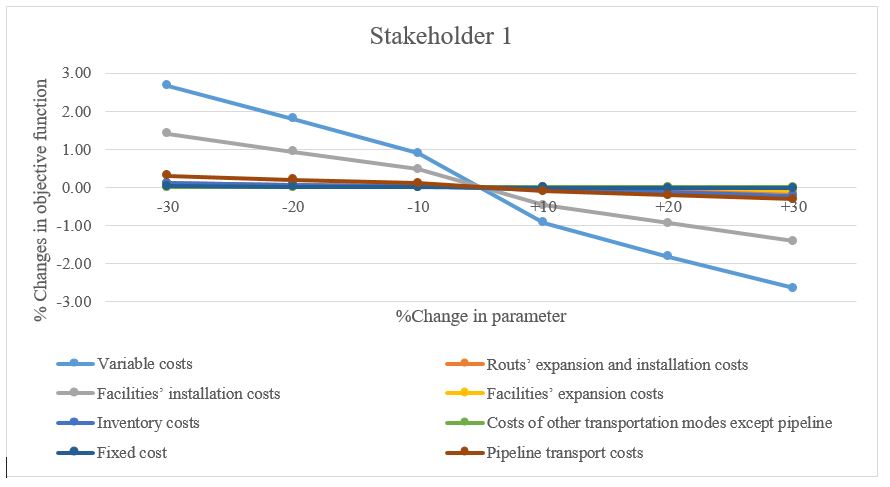}
    \caption{Rate of changes in objective function values of the first stakeholder}
    \label{fig:fig3}
\end{figure}
\begin{figure}[h]
    \centering
    \includegraphics[width=1\linewidth,trim={0.2in 0.2in 0.2in 0.2in},clip]{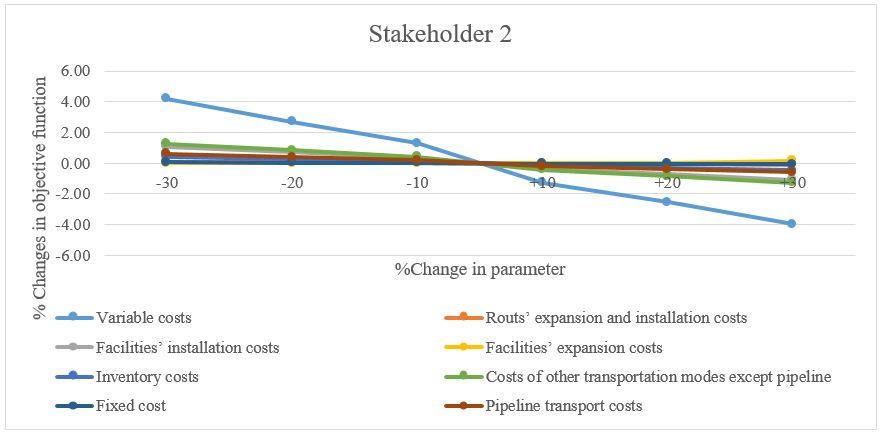}
    \caption{Rate of change in objective function values of the second stakeholder}
    \label{fig:fig4}
\end{figure}
\begin{figure}[h]
    \centering
    \includegraphics[width=1\linewidth,trim={0.2in 0.2in 0.2in 0.2in},clip]{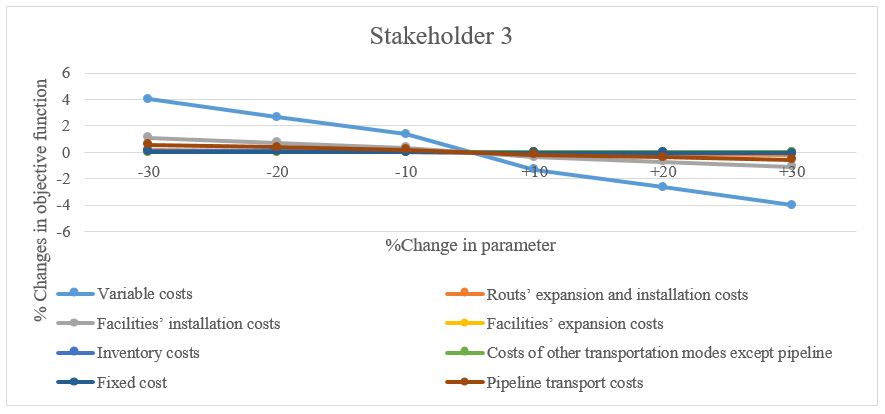}
    \caption{Rate of change in objective function values of the third stakeholder}
    \label{fig:fig5}
\end{figure}

The percentage changes of the objective function for each of the first, second and third stakeholders are shown in Figures \ref{fig:fig3}, \ref{fig:fig4} and \ref{fig:fig5}, respectively.

As seen, with a 30\% change in parameters, for all three stakeholders, the highest change in the objective function value is in the 2.5-4\% range, which corresponds to the variable cost group. In addition to variable costs, the highest changes in the objective function correspond to the variations in facility construction and transportation costs. The reason for the higher change in the objective function with changes in variable costs is the associated construction and transportation costs, which make up the highest share of costs compared with other costs, according to Figure \ref{fig:fig1}.

\section{Conclusions}\label{sec:conclusions}

Multi- period, multi objective, multi-level, multi-stakeholder and multi-product DSCPSC model has been investigated considering the competition between different stakeholders using job creation maximization, profit maximization and pollution minimization objective functions. The proposed model is an MILP model. In addition to different stakeholders, the dynamic design of the supply chain has also been considered. Dynamic design of the supply chain simultaneously considers all three approaches of new facility construction, expansion of existing facilities and capacity reduction and/or facility close down. The other aspect of dynamism is related to the distribution centers in this model. Distribution centers are one of the sources shared between stakeholders. Different stakeholders can use a distribution center within a period of time and pay its expenses according to their use of the distribution center. In addition to the optimal use of the capacity of a distribution center, another advantage of this approach is the minimization of each stakeholder’s expenses. Considering that the model is a multi-objective one, fuzzy approach has been used to convert it to a single-objective model. The proposed model was executed using GAMS 24.1.2 software for a part of the real PSC data. Sensitivity analysis was then performed to determine the sensitivity of the results to the variation of the coefficients of objective function (cost parameters). According to the sensitivity analysis, the highest variations of the objective function were observed with respect to the changes in coefficients of variable, facility construction and transportation costs. The reason for this observation is that these coefficients make up a great part of the whole costs for each of the stakeholders.












\bibliographystyle{IEEEtran}
\bibliography{main.bib}

\section*{Appendix A: parameters, sets and variables used in the mode}\label{sec:AppendixA}

\begin{table}[H]
   
\captionsetup{labelformat=empty}
\centering
\caption*{A-1: Sets}
\label{A1}
\begin{tabular}{|l |l |}

\hline

$K$& Set of existing refineries	 \tabularnewline\hline

$\acute{K}$&Set of new refineries\tabularnewline\hline	 

$L$&Set of existing DCs\tabularnewline\hline	 

$\acute{L}$&Set of the new DCs	 \tabularnewline\hline

$M$&Set of costumer zone\tabularnewline\hline	 

$V$&Set of transportation modes	 

\tabularnewline\hline

$LCV$&Set of capacity of transportation modes\tabularnewline\hline	 
$P$&Set of products	\tabularnewline\hline 

$EK$&Set of capacity of new refineries	 \tabularnewline\hline

$EL$&Set of capacity of new DCs\tabularnewline\hline	 

$UK$&Set of capacity expansion of existing refineries	\tabularnewline\hline 

$UL$&Set of capacity expansion of existing DCs	 \tabularnewline\hline

$EV$&Set of capacity expansion of existing pipelines	\tabularnewline\hline 

$EZ$&Set of capacity of storage tanks in DC	 \tabularnewline\hline

$LV$&Set of capacity of new pipeline\tabularnewline\hline	 

$RV$&Set of new pipeline routes 	 \tabularnewline\hline

$E$&Set of stakeholders including government and private sectors	 \tabularnewline\hline

$EN$&Set of regions	EN\tabularnewline\hline
$LEV$&Set of education levels	LEV\tabularnewline\hline
$T$&Set of Time of periods	 

\tabularnewline\hline

\end{tabular}
\end{table}

\begin{table}[H]
   
\captionsetup{labelformat=empty}
\centering
\caption*{A-2: Indices}
\label{A2}
\begin{tabular}{|l |l |}

\hline

$k$&Refinery\tabularnewline\hline
$l$&DC\tabularnewline\hline
$m$&Costumer zone\tabularnewline\hline
$v$&Transportation mode\tabularnewline\hline
$lcv$&Capacity of transportation mode\tabularnewline\hline
$p$&Product\tabularnewline\hline
$ek$&Capacity of new refinery\tabularnewline\hline
$el$&Capacity of new DC\tabularnewline\hline
$uk$&Capacity expansion of existing refinery\tabularnewline\hline
$ul$&Capacity expansion of existing DC\tabularnewline\hline
$ev$&Capacity expansion of existing pipeline\tabularnewline\hline
$ez$&Capacity of the storage tank in DC\tabularnewline\hline
$lv$&Capacity of new pipeline\tabularnewline\hline
$rv$&New pipeline route\tabularnewline\hline
$e$&Beneficiaries including government and private sectors\tabularnewline\hline
$en$&regions\tabularnewline\hline
$lev$&Education levels\tabularnewline\hline
$PS$&Private sectors\tabularnewline\hline
$G$&Government\tabularnewline\hline
$t$&Time of periods

\tabularnewline\hline

\end{tabular}
\end{table}

\begin{table}[H]
   
\captionsetup{labelformat=empty}
\centering
\caption*{A-3: Parameters}
\label{A3}
\begin{tabular}{|l |p{6.5cm} |}

\hline
 
  $d_{p}^{t,e}$&Demand for product p within time period t supplied by stakeholder e ($e \in E,t\in T,p\in P$)
 \tabularnewline\hline

$d_{p,m}^{t}$&Demand of customer m for product p within time period t ($m \in M,t\in T,p\in P$)
 \tabularnewline\hline

$capl_{l}^{ul}$&Capacity expansion for existing distribution center 1 in level ul ($l\in L, ul\in L$)
\tabularnewline\hline
 $capk_{k}^{uk}$&Capacity expansion for existing refinery k in level uk ($k \in K,uk\in UK$)
 
 \tabularnewline\hline
 $ick_{k}^{}$&Initial capacity of existing refinery k   ($k \in K$)
 \tabularnewline\hline
 $icl_{p,l}^{}$&Initial capacity of existing distribution center 1 ($l \in L,p\in P$)
  \tabularnewline\hline

\end{tabular}
\end{table}

\begin{table}[H]
   
\captionsetup{labelformat=empty}
\centering
\caption*{A-3: Parameters (Continued)}
\label{A3-1}
\begin{tabular}{|l |p{6.5cm} |}

\hline

$Mk_{\acute{k}}^{ek}$&Minimum rate of capacity utilization coefficient of new refinery k in level ek ($ek \in EK,k\in \acute{K}$) \tabularnewline\hline

$Ml_{p,\acute{l}}^{ez}$&Minimum rate of capacity utilization coefficient of new distribution center l for product p in level el ($ez \in EZ,p \in P,l\in \acute{L}$) \tabularnewline\hline

$ivk_{k}^{0}$&Initial inventory level in refinery k ($k\in (\acute{K} \cup K)$)\tabularnewline\hline

$ivl_{l}^{0}$&Initial inventory level in distribution center l ($l\in (\acute{L} \cup L)$)\tabularnewline\hline

$\lambda E_{en}^{}$&Coefficient of pollution significance in area en ($en\in EN$) \tabularnewline\hline

$Nct_{\acute{l}}^{ez}$&Capacity of storage tank in level ez in new distribution center l ($l\in \acute{L},ez \in EZ$)\tabularnewline\hline

$Ncl_{l}^{el}$&Capacity of new distribution center l in level el ($l\in \acute{L},el \in EL$)\tabularnewline\hline

$Nck_{\acute{k}}^{ekl}$&Capacity of new refinery k in level ek   ($k\in \acute{K},ek \in EK$)\tabularnewline\hline

$icapkl_{k,l}^{}$&Initial capacity of the pipeline between refinery k and distribution center ($k\in K,l \in L$)\tabularnewline\hline

$icaplpl_{lp,l}^{}$&Initial capacity of the pipeline between distribution center l and distribution center ($lp,l\in L,lp \neq l$)\tabularnewline\hline

$trc_{v}^{lcv}$&Capacity of transportation mode v in level lcv  ($v\in V,lcv\neq LCV$)\tabularnewline\hline

$capkl_{k,l}^{ev}$&Capacity expansion of the route between existing refinery k and distribution center l in level ev ($k\in K,l \in L, ev \in EV$)\tabularnewline\hline

$caplpl_{lp,l}^{ev}$&Capacity expansion of the route between existing distribution center lp and existing distribution center l in level ev ($lp,l\in L,lp \neq l, ev \in EV$)\tabularnewline\hline

$TPP$&Time period\tabularnewline\hline

$clv_{lv}^{}$&Capacity of pipeline transportation mode in level lv ($lv \in LV$)\tabularnewline\hline

$Rkl_{k,l}^{}$&Zero and one matrix indicating the existing routes between existing refinery k and existing distribution center l ($l\in L,k \in K$)\tabularnewline\hline
 
$Rlpl_{lp,l}^{}$&Zero and one matrix indicating the existing routes between existing distribution center lp and existing distribution center l ($l,lp\in L,lp \neq l$)\tabularnewline\hline
 
$nmax_{v}^{lcv,t}$&Maximum number of available transportation mode v in level lcv in time t ($lcv\in LCV,v \in V, t \in T$)\tabularnewline\hline

$w_{i}^{}$&Maximum crude oil production from oil field i   ($i \in I$)\tabularnewline\hline

$\mu_{p}^{}$&Production rate of product p from one crude oil barrel ($p \in P$)\tabularnewline\hline

$lk_{k}^{}$&Minimum inventory level of refinery k ($k \in (K \cup \acute{K}) P$)\tabularnewline\hline
   
$ll_{l}^{}$&Minimum inventory level of distribution center l ($l \in (L \cup \acute{L}) P$)\tabularnewline\hline

$M$&A large number \tabularnewline\hline
   
$xcostk_{k}^{ek,t}$&Installation cost of new refinery k in level ek during time t ($t \in T,ek \in EK, k \in \acute{K}) P$)\tabularnewline\hline

$xcostl_{l}^{el,t}$&Installation cost of new distribution center l in level el during time t  ($l \in  \acute{L}, t \in T$)\tabularnewline\hline

$ucostk_{k}^{uk,t}$&Expansion cost of existing refinery k in level uk during time t ($t \in T, p \in P, k \in K,uk \in UK) P$)\tabularnewline\hline

$ucostl_{l}^{ul,p,t}$&Expansion cost of existing distribution center l in level ul for product p during time t ($t \in T, p \in P, l \in L,ul \in L) P$)\tabularnewline\hline

$ycostkl_{k,l}^{ev,t}$&Expansion cost of pipeline between existing refinery k and existing distribution center l in level lv during time t ($t \in T, ev \in EV, k \in K, l \in L) P$)\tabularnewline\hline

$ycostlpl_{lp,l}^{ev,t}$&Expansion cost of pipeline between existing distribution center lp and existing distribution center l in level lv during time t ($t \in T, ev \in EV, lp,l \in L, l \neq lp$)\tabularnewline\hline

$hcostk_{k}^{}$&Storage cost of crude oil in refinery k during time t ($t \in T, k \in (K \cup \acute{K})$)\tabularnewline\hline

$hcostl_{p,l}^{}$&Storage cost of product p in distribution center l during time t ($t \in T,p \in P, l \in (L \cup \acute{L})$)\tabularnewline\hline
$qcostkl_{k,l}^{t}$&Flow cost of crude oil between refinery k and distribution center l during time t ($l \in (L \cup \acute{L}),k \in (K \cup \acute{K}), t \in T$)\tabularnewline\hline
$qcostlpl_{lp,l}^{t}$&Flow cost of products between existing distribution center lp and existing distribution center l in level lv during time ($t \in T, lp \in L, l \in (L \cup \acute{L}) OR  t \in T, lp \in \acute{L}, l \in \acute{L}$)\tabularnewline\hline

\end{tabular}
\end{table}

\begin{table}[H]
   
\captionsetup{labelformat=empty}
\centering
\caption*{A-3: Parameters (Continued)}
\label{A3-2}
\begin{tabular}{|l |p{6cm} |}

\hline

$rcostkl_{k,l,v}^{lv,rv,t}$&Installation cost of new pipeline rv between refinery k and distribution center l in level lv during time t ($lv \in LV, rv \in RV, l \in (L \cup \acute{L}) OR t \in T, lp \in \acute{L}, l \in \acute{L}$)\tabularnewline\hline

$ycostkl_{k,l}^{ev,t}$&Expansion cost of pipeline between existing refinery k and existing distribution center l in level lv during time t ($t \in T, ev \in EV, k \in K, l \in L $)\tabularnewline\hline

$ycostlpl_{lp,l}^{ev,t}$&Expansion cost of pipeline between existing distribution center lp and existing distribution center l in level lv during time t ($t \in T, ev \in EV, lp,l \in L, l \neq lp$)\tabularnewline\hline

$hcostk_{k}^{}$&Storage cost of crude oil in refinery k during time t ($t \in T, k \in (K \cup \acute{K})$)\tabularnewline\hline

$hcostl_{p,l}^{}$&Storage cost of product p in distribution center l during time t ($t \in T,p \in P, l \in (L \cup \acute{L})$)\tabularnewline\hline

$qcostkl_{k,l}^{t}$&Flow cost of crude oil between refinery k and distribution center l during time t ($t \in T, l \in (L \cup \acute{L}), k \in (K \cup \acute{K})$)\tabularnewline\hline

$qcostlpl_{lp,l}^{t}$&Flow cost of products between existing distribution center lp and existing distribution center l in level lv during time ($t \in T,lp \in L, l \in (L \cup \acute{L}) OR\;t \in T,lp \in \acute{L}, l \in \acute{L}$)\tabularnewline\hline

$rcostkl_{k,l,v}^{lv,rv,t}$&Installation cost of new pipeline rv between refinery k and distribution center l in level lv during time t ($lv \in LV, rv \in RV,t \in T, l \in (L \cup \acute{L}), k \in (K \cup \acute{K})$v=pipeline mode)\tabularnewline\hline

$rcostlpl_{lp,l,p}^{lv,rv,t}$&Installation cost of new pipeline rv between other distribution centers lp and the corresponding distribution center l in level lv during time t  ($lv \in LV, rv \in RV,t \in T, l,lp \in (L \cup \acute{L})$)\tabularnewline\hline

$ncostkl_{k,l,v}^{lcv,t}$&Unit cost of transportation mode v between refinery k and distribution center l in capacity level lcv during time t($lcv \in LCV, v \in V,t \in T, l \in (L \cup \acute{L}), k \in (K \cup \acute{K})$)\tabularnewline\hline

$ncostlm_{l,m,v}^{lcv,t}$&customer area m in capacity level lcv during time t  ($lcv \in LCV, v \in V,t \in T, l \in (L \cup \acute{L}), m\in M$)\tabularnewline\hline

$ncostlpl_{lp,l,v}^{lcv,t}$&Unit cost of transportation mode v between distribution centers l and lp in capacity level lcv during time t ($lcv \in LCV, v \in V,t \in T, l,lp \in (L \cup \acute{L})$)\tabularnewline\hline

$icost_{p}^{t}$&Importation cost of product p during time t   ($t \in T,p \in P$)\tabularnewline\hline

$ncostl_{l}^{ez,t}$&Installation cost of storage tank with capacity ez in distribution center l during time t ($t \in T,ez \in EZ,l \in \acute{L}$)\tabularnewline\hline

$pcostk_{k}^{t}$&Operation cost of refinery k during time t ($t \in T,k \in (K \cup \acute{K})$)\tabularnewline\hline

$pcostl_{p,l}^{t}$&Operation cost of product p in distribution center l during time t ($t \in T,l \in (L \cup \acute{L}),p \in P$)\tabularnewline\hline

$Fcostk_{k}^{t}$&Fixed operation cost of refinery k during time t ($t \in T,k \in (K \cup \acute{K})$)\tabularnewline\hline

$Fcostl_{l}^{t}$&Fixed operation cost of distribution center l during time t ($t \in T,l \in (L \cup \acute{L})$)\tabularnewline\hline

$OP_{}^{t}$&Crude oil price during time t($t \in T$)\tabularnewline\hline

$ERPP_{p}^{t,e}$&Exportation cost of product p during time t by stakeholder e ($t \in T,p \in P$)\tabularnewline\hline

$WCost_{lev}^{t}$&Labor force cost with education level lev during time t ($t \in T$)\tabularnewline\hline

$clcostk_{k}^{t}$&Closing down cost of refinery during time t ($t \in T,k \in K$)\tabularnewline\hline

$Nk_{en}^{k}$&Zero and one matrix indicating area en in which candidate refinery k has been constructed  ($en \in EN,k \in \acute{K}$)\tabularnewline\hline

$Nek_{en}^{k}$&Zero and one matrix indicating area en in which existing refinery k has been expanded  ($en \in EN,k \in K$)\tabularnewline\hline

$Nl_{en}^{l}$&Zero and one matrix indicating area en where candidate distribution center l has been constructed ($en \in EN,l \in \acute{ L}$)\tabularnewline\hline

$Nel_{en}^{l}$&Zero and one matrix indicating area en, in which existing distribution center l has been constructed  ($en \in EN,l \in \acute{ L}$)\tabularnewline\hline

\end{tabular}
\end{table}

\begin{table}[H]
   
\captionsetup{labelformat=empty}
\centering
\caption*{A-3: Parameters (Continued)}
\label{A3-2}
\begin{tabular}{|l |p{6cm} |}

\hline

$Pulk$&Pollution per barrel of crude oil refined\tabularnewline\hline

$Pull$&Pollution in distribution centers\tabularnewline\hline

$Pulv_{lcv}^{v}$&Pollution per kilometer of transportation mode v with capacity level lcv ($lcv \in LCV,v \in V$)\tabularnewline\hline

$Per$&Pollution coefficients for capacity expansion of refinery and distribution center\tabularnewline\hline

$dis_{m}^{l}$&Distance between distribution center l and customer area m ($m \in M,l \in (L \cup \acute{L}) $)\tabularnewline\hline

$W_{en}^{\acute{en}}$&Coefficient of significance of labor force transfer from area en to area ($en \in EN,en \in \acute{EN}$)\tabularnewline\hline

$Nlab_{en}^{lev}$&New labor force with education level lev in area en ($en \in EN,lev \in LEV$)\tabularnewline\hline

$WNK_{lev}^{}$&Required labor force with education level lev for a new refinery($lev \in LEV$)\tabularnewline\hline

$WEK_{lev}^{}$&Required labor force with education level lev for capacity expansion of existing refinery ($lev \in LEV$)\tabularnewline\hline

$WNL_{lev}^{}$&Required labor force with education level lev for a new distribution center ($lev \in LEV$)\tabularnewline\hline

$WEL_{lev}^{}$&Required labor force with education level lev for capacity expansion of existing distribution center ($lev \in LEV$)\tabularnewline\hline

$Maxnum_{en}^{}$&Maximum construction and expansion of refineries and distribution centers in each area ($en \in EN$)\tabularnewline\hline

\end{tabular}
\end{table}

\begin{table}[H]
   
\captionsetup{labelformat=empty}
\centering
\caption*{A-5: Positive integer variables}
\label{A5}
\begin{tabular}{|l |p{6cm} |}

\hline

$nkl_{p,k,l,v}^{lcv,t,e}$& Number of fleet of transportation mode v at capacity level lcv between refinery k and DC l and by stakeholder e during t ($t \in T,lcv \in LCV, e \in E,v \in V, k \in (K \cup \acute{K}), l \in (L \cup \acute{L})$)\tabularnewline\hline

$nlpl_{p,lp,l,v}^{lcv,t,e}$& Number of fleet of transportation mode v at capacity level lcv between DC lp and DC l by stakeholder e during t ($t \in T,lcv \in LCV, e \in E,v \in V, p \in P, l,lp \in (L \cup \acute{L})$)\tabularnewline\hline

$nlm_{p,l,m,v}^{lcv,t,e}$& Number of fleet of transportation mode v at capacity level lcv between DC l and customer zone m by stakeholder e during t ($t \in T,lcv \in LCV, e \in E,v \in V, p \in P, l,lp \in (L \cup \acute{L})$)\tabularnewline\hline

$n_{p,l}^{ez,t,e}$& The number of storage tanks in capacity level ez in product p in DC l in by stakeholder e during t ($t \in T, e \in E, ez \in EZ, p \in P, l \in (L \cup \acute{L})$)\tabularnewline\hline

$HENK_{en,k}^{\acute{en},lev,t}$& The number of labor forces in level lev, which are worked in  new refinery k in region en from region   in capacity level lcv during t ($lev \in LEV,t \in T, en,\acute{en} \in EN, k \in \acute{K}$)\tabularnewline\hline

$HEEK_{en,k}^{\acute{en},lev,t}$& TThe number of labor forces in level lev, which are worked in refinery k in region en   from region   in capacity level lcv during t ($t \in T, en,\acute{en} \in EN, lev \in LEV, k \in K$)\tabularnewline\hline

$HENl_{en,l}^{\acute{en},lev,t}$&The number of labor forces in level lev, which are worked in  DC l in region en from region   in capacity level lcv during t ($t \in T, en,\acute{en} \in EN, lev \in LEV, l \in \acute{L}$)\tabularnewline\hline

$HENl_{en,l}^{\acute{en},lev,t}$&The number of labor forces in level lev, which are worked in  DC l in region en from region   in capacity level lcv during t ($t \in T, en,\acute{en} \in EN, lev \in LEV, l \in L$)\tabularnewline\hline

\end{tabular}
\end{table}

\begin{table}[H]
   
\captionsetup{labelformat=empty}
\centering
\caption*{A-4: Positive variables}
\label{A4-0}
\begin{tabular}{|l |p{6cm} |}

\hline

$qkl_{p,k,l,v}^{t,e}$&Rate of crude oil stream from refinery k to distribution center l via transportation mode v by stakeholder e during time t ($t \in T,p \in P, e \in E,v \in V, k \in (K \cup \acute{K}), l \in (L \cup \acute{L})$)\tabularnewline\hline

$qlm_{p,l,m,v}^{t,e}$&Rate of product stream from distribution center l to customer m via transportation mode v by stakeholder e during time t ($t \in T,p \in P, e \in E,v \in V, k \in (K \cup \acute{K}), m \in M$)\tabularnewline\hline

\end{tabular}
\end{table}

\begin{table}[H]
   
\captionsetup{labelformat=empty}
\centering
\caption*{A-4: Positive variables}
\label{A4}
\begin{tabular}{|l |p{6cm} |}

\hline

$qlpl_{p,lp,l,v}^{t,e}$&Rate of product stream from other distribution centers lp to distribution center l via transportation mode v by stakeholder e during time t ($t \in T, e \in E,v \in V,l,lp \in (L \cup \acute{L})$)\tabularnewline\hline

$i_{p,l}^{t,e}$&Imports of product p from distribution center l by stakeholder e during time t ($t \in T,p \in P, e \in E, l \in (L \cup \acute{L}) $)\tabularnewline\hline

$vl_{p,l}^{t,e}$&Inventory of product p in distribution center l by stakeholder e during time t ($t \in T,p \in P, e \in E, l \in (L \cup \acute{L}) $)\tabularnewline\hline

$vk_{k}^{t,e}$& Inventory in refinery k by stakeholder e during time t ($t \in T, e \in E, k \in (K \cup \acute{K}) $)\tabularnewline\hline

$Ep_{p,l}^{t,e}$& Exports of product p from distribution center l by stakeholder e during time t ($t \in T,p \in P, e \in E, l \in (L \cup \acute{L}) $)\tabularnewline\hline

$\xi_{l}^{t,e}$&Percentage of the usage of distribution center l by stakeholder e during time t ($t \in T, e \in E, l \in L$)\tabularnewline\hline

$RLab_{en}^{lev,t}$& Number or required labor forces in level lev during time t($t \in T, en \in EN, lev \in LEV$)\tabularnewline\hline

$ALab_{en}^{lev,t}$& Available labor force with education level lev in area en during time t ($t \in T, en \in EN, lev \in LEV$)\tabularnewline\hline

\end{tabular}
\end{table}

\begin{table}[H]
   
\captionsetup{labelformat=empty}
\centering
\caption*{A-6: Binary variables}
\label{A6}
\begin{tabular}{|p{8.3cm} |}

\hline

  \[\label{cons93}
    xk_{k}^{ek,t,e}=\left\{
                \begin{array}{ll}
                  1\;\;\;\;\;If \;the\; refinery\; k\; is\; installed \;in \;level\; ek \\\;\;\;\;\;\;\;during\;time \;t \;by \;e\\
                  0\;\;\;\;\;Otherwise\;\;\;\;(t \in T, ek \in EK, e \in E,k \in \acute{K})\nonumber
                \end{array}
              \right.
  \]

\tabularnewline\hline

  \[\label{cons93}
    \tau l_{p,l}^{ul,t,e}=\left\{
                \begin{array}{ll}
                  1\;\;\;\;\;If \;the\; distribution\;center\; l\; is\; expanded \; \\\;\;\;\;\;\;in \;level\; ul\;during\;time \;t \;by \;e\\
                  0\;\;\;\;\;Otherwise(t \in T, ul \in UL, e \in E,l \in L,p \in P)\nonumber
                \end{array}
              \right.
  \]

\tabularnewline\hline

  \[\label{cons93}
    \tau k_{k}^{uk,t,e}=\left\{
                \begin{array}{ll}
                  1\;\;\;\;\;If \;the\; refinery\;k\; is\; expanded \;in \;level\;uk \\\;\;\;\;\;\; \;during\;time \;t \;by \;e\\
                  0\;\;\;\;\;Otherwise\;\;\;\;(t \in T, uk \in UK, e \in E,k \in K)\nonumber
                \end{array}
              \right.
  \]

\tabularnewline\hline

  \[\label{cons93}
    ykl_{k,l}^{ev,t,e}=\left\{
                \begin{array}{ll}
                  1\;\;\;\;\;If \;existing\; route\;between\; existing\; refinery  \\\;\;\;\;\;\; k \;and\;DC\;is\; expanded\; with\; level\; ev \;in \\\;\;\;\;\;\;period\; t  \;by \;e\\
                  0\;\;\;\;\;Otherwise(t \in T, ev \in EV, e \in E,k \in K,l \in L)\nonumber
                \end{array}
              \right.
  \]

\tabularnewline\hline

  \[\label{cons93}
    ylpl_{lp,l}^{ev,t,e}=\left\{
                \begin{array}{ll}
                  1\;\;\;\;\;If \;existing\; route\;between\; existing\; DC  \\\;\;\;\;\;\; lp \;and\;existing\;DC\;is\; expanded\; with\; level\;  \\\;\;\;\;\;\;ev \;in\;period\; t  \;by \;e\\
                  0\;\;\;\;\;Otherwise(t \in T, ev \in EV, e \in E,l,lp \in L)\nonumber
                \end{array}
              \right.
  \]
  
  \tabularnewline\hline
    \[\label{cons93}
    z_{l}^{ez,t,e}=\left\{
                \begin{array}{ll}
                  1\;\;\;\;\;If \;the\; number\;of\; storage\; tanks\; is \;fixed  \\\;\;\;\;\;\;with\;level\;ez\; in\; new\; level\;  \\\;\;\;\;\;\;DC \;l\;in\;period\; t  \;by \;e\\
                  0\;\;\;\;\;Otherwise(t \in T, ez \in EZ, e \in E,\acute{l} \in L)\nonumber
                \end{array}
              \right.
  \]

\tabularnewline\hline
    \[\label{cons93}
    rkl_{k,l,v}^{lv,rv,t,e}=\left\{
                \begin{array}{ll}
                  1\;\;\;\;\;If \;route\; rv\;between\; refinery\; k \; and \;DC\;l\; \\\;\;\;\;\;\;is\;installed\; with\; level\; lv\; \;in \;period\; t\;by \;e \\
              0\;\;\;\;\;Otherwise\\\;\;\;\;(t \in T, ev \in EV, k \in (K \in \acute{K}),l \in \acute{L} \in L)\nonumber
                \end{array}
              \right.
  \]
\tabularnewline\hline

  \[\label{cons93}
    rlpl_{lp,l,v}^{lv,rv,t,e}=\left\{
                \begin{array}{ll}
                  1\;\;\;\;\;If \;route\; rv\;between\; DC\; lp \; and \;DC\;l\; \\\;\;\;\;\;\;is\;installed\; with\; level\; lv\; \;in \;period\; t\;by \;e \\
              0\;\;\;\;\;Otherwise(t \in T, ev \in EV,l,lp \in (\acute{L} \in L)\nonumber
                \end{array}
              \right.
  \]
\tabularnewline\hline

    \[\label{cons93}
    rlp\acute{l}_{lp,l,v}^{t}=\left\{
                \begin{array}{ll}
                  1\;\;\;\;\;If \;route\; rv\;between\; DC\; lp \; and \;DC\;l\; \\\;\;\;\;\;\;is\;installed\; in \;period\; t\;by \;e \\
              0\;\;\;\;\;Otherwise(t \in T, e \in E,l,lp \in (\acute{L} \in L)\nonumber
                \end{array}
              \right.
  \]
\tabularnewline\hline

\end{tabular}
\end{table}

\begin{table}[H]
   
\captionsetup{labelformat=empty}
\centering
\caption*{A-6: Binary variables (Continued)}
\label{A6-1}
\begin{tabular}{|p{8.3cm} |}

\hline

    \[\label{cons93}
    S_{lp,l}^{t}(Auxiliary\;variable)=\left\{
                \begin{array}{ll}
                  1\;\;\;\;\;If \;there\; is\;a\; flow\; from \; DC  \\\;\;\;\;\;\;\;lp\;to\;DC\;l\; in \;period\; t\\
              0\;\;\;\;\;Otherwise(l,lp \in (\acute{L} \in L)\nonumber
                \end{array}
              \right.
  \]
\tabularnewline\hline

    \[\label{cons93}
    \psi k_{k}^{t,e}=\left\{
                \begin{array}{ll}
                  1\;\;\;\;\;If \;k\; is\;closed\; in\; period \; t\;\;\;\;\;\;\;\;\;\; \;\;\;\;\;\;\;\;\;\;\;\;\;\;\;\;\;\;\;\;\;\;\;\;\\
              0\;\;\;\;\;Otherwise(t \in T,k \in K,e \in E)\nonumber
                \end{array}
              \right.
  \]
\tabularnewline\hline

\end{tabular}
\end{table}

\label{sec:AppendixB}

\begin{table*}[htbp]
\captionsetup{labelformat=empty}
\centering
\caption*{Appendix B: Optimal values in the case study}   

\label{table2}
\centering
\begin{tabular}{|>{\centering}c |>{\centering}c |>{\centering}c |>{\centering}c |>{\centering}c |}

\hline
 \multicolumn{5}{|c|}{$xk_{k}^{ek,t}, k \in \acute{k}$\;\;\;\; (Binary variable)}\\
\hline
   ek  &   k &  t &  e &  Value \tabularnewline \hline
   1&	3&	1&	3&	1\tabularnewline \hline
   1&	2&	1&	2&	1\tabularnewline \hline
   1&	1&	1&	1&	1\tabularnewline

\hline

\end{tabular}
\end{table*}


\begin{table*}[htbp]
   
\captionsetup{labelformat=empty}
\label{resulttable}
\centering
\begin{tabular}{|>{\centering}c |>{\centering}c |>{\centering}c |>{\centering}c|>{\centering}c|>{\centering}c |>{\centering}c|>{\centering}c |>{\centering}c|>{\centering}c|>{\centering}c |>{\centering}c |}

\hline
 \multicolumn{12}{|c|}{$nsu$}\\
\hline
t&	p&	value&	t&	p&	value&	t&	p&	value&	t&	p&	value\tabularnewline \hline
1&	1&	127.83&	1&	2&	48.748&	2&	1&	151.88&	2&	2&	60.114\tabularnewline \hline
3&	1&	183.33&	3&	2&	63.028&	&	&		&	&	&\tabularnewline

\hline

\end{tabular}
\end{table*}


\begin{table*}[htbp]
   
\captionsetup{labelformat=empty}

\label{resulttable}
\centering
\begin{tabular}{|>{\centering}c |>{\centering}c |>{\centering}c |>{\centering}c|>{\centering}c|>{\centering}c |>{\centering}c|>{\centering}c |>{\centering}c|>{\centering}c|>{\centering}c |>{\centering}c |}

\hline
 \multicolumn{12}{|c|}{$\psi_{k}^{t,e}$\;\;\;\; }\\
\hline

t&	e&	k&	value&	t&	e&	k&	value&	t&	e&	k&	value\tabularnewline \hline
3&	3&	3&	1&	    3&	3&	2&	   1&	3&	3&	1&	1\tabularnewline \hline
2&	2&	2&	1&	    2&	1&	3&	1&	2&	1&	1&	1\tabularnewline \hline
1&	1&	3&	1&   	1&	1&	2&	1&	1&	1&	1&	1\tabularnewline \hline

\end{tabular}
\end{table*}


\begin{table*}[htbp]
   
\captionsetup{labelformat=empty}
\label{resulttable}
\centering
\begin{tabular}{|>{\centering}c |>{\centering}c |>{\centering}c |>{\centering}c|>{\centering}c|>{\centering}c|>{\centering}c|>{\centering}c |}

\hline
 \multicolumn{8}{|c|}{$rkl_{k,l,v}^{lv,rv,t}, k \in K, l \in L$\;\;\;\; (Binary variable)}\\
\hline
   k&	l&	lv&	rv&	v&	t&	e&	value\tabularnewline \hline
   3&	3&	8&	1&	3&	1&	3&	1\tabularnewline \hline
   2&	3&	8&	1&	3&	1&	2&	1\tabularnewline \hline
   2&	2&	8&	1&	3&	1&	2&	1\tabularnewline \hline
   1&	2&	8&	1&	3&	1&	1&	1\tabularnewline \hline
   1&	1&	8&	1&	3&	1&	1&	1\tabularnewline

\hline

\end{tabular}
\end{table*}



\begin{table*}[hbt!]
   
\captionsetup{labelformat=empty}
\label{resulttable}
\centering
\begin{tabular}{|>{\centering}c |>{\centering}c |>{\centering}c |>{\centering}c|>{\centering}c|>{\centering}c |>{\centering}c|>{\centering}c |>{\centering}c|>{\centering}c|>{\centering}c |>{\centering}c |>{\centering}c |>{\centering}c|>{\centering}c|>{\centering}c |>{\centering}c |>{\centering}c |>{\centering}c|>{\centering}c|}

\hline
 \multicolumn{20}{|c|}{$vl_{p,l}^{t}, l \in L$\;\;\;\; (1000 barrel\;per\;day)}\\
\hline
p&	l&	t&	e&	value&	p&	l&	t&e&	value&	p&	l&	t&	e&	value&	p&	l&	t&	e&	value\tabularnewline \hline
2&	2&	1&	1&	35.192&	2&	2&	2&1&	28.559&	1&	1&	2&	1&	21.98&	1&	2&	1&	2&	19.294\tabularnewline \hline
2&	3&	2&	3&	4.9	&	&	&	&&		&	&	&	&	&		&	&	&    &   &\tabularnewline

\hline

\end{tabular}
\end{table*}


\begin{table*}[hbt!]
   
\captionsetup{labelformat=empty}
\label{resulttable}
\centering
\begin{tabular}{|>{\centering}c |>{\centering}c |>{\centering}c |>{\centering}c|>{\centering}c|>{\centering}c |>{\centering}c|>{\centering}c |>{\centering}c|>{\centering}c|>{\centering}c |>{\centering}c |>{\centering}c |>{\centering}c|>{\centering}c|>{\centering}c |>{\centering}c |>{\centering}c |>{\centering}c|>{\centering}c|>{\centering}c|}

\hline
 \multicolumn{21}{|c|}{$qkl_{p,k,l,v}^{t}, l \in \acute{L}, k \in \acute{K}$\;\;\;\; (1000 barrel\;per\;day)}\\
\hline
p&	k&	l&	v&	t&	e&	value&	p&	k&	l&	v&	t&	e&	value&	p&	k&	l&	v&	t&	e&	value\tabularnewline \hline
1&	3&	3&	3&	3&	3&	65.407&	1&	3&	3&	3&	1&	3&	65.407&	1&	2&	3&	3&	2&	2&	49.1\tabularnewline \hline
1&	2&	2&	3&	3&	2&	47.997&	1&	1&	1&	3&	1&	1&	43.605&	1&	1&	2&	3&	3&	1&	34.55\tabularnewline \hline
2&	3&	3&	3&	2&	3&	22.363&	1&	2&	3&	3&	3&	2&	17.41&	2&	2&	2&	3&	2&	2&	5.64\tabularnewline \hline

\end{tabular}
\end{table*}


\begin{table*}[hbt!]
   
\captionsetup{labelformat=empty}
\label{resulttable}
\centering
\begin{tabular}{|>{\centering}c |>{\centering}c |>{\centering}c |>{\centering}c|>{\centering}c|>{\centering}c |>{\centering}c|>{\centering}c |>{\centering}c|>{\centering}c|>{\centering}c |>{\centering}c |>{\centering}c |>{\centering}c|>{\centering}c|>{\centering}c |>{\centering}c |>{\centering}c |>{\centering}c|>{\centering}c|>{\centering}c|}

\hline
 \multicolumn{21}{|c|}{$qkl_{p,k,l,v}^{t}, l \in L, k \in \acute{K}$\;\;\;\; (1000 barrel\;per\;day)}\\
\hline

p&	k&	l&	v&	t&	e&	value&	p&	k&	l&	v&	t&	e&	value&	p&	k&	l&	v&	t&	e&	value\tabularnewline \hline
1&	1&	1&	3&	2&	1&	200	&   1&	1&	1&	3&	3&	1&	158.0&	1&	3&	3&	3&	2&	3&	150.9\tabularnewline \hline
1&	2&	2&	3&	1&	2&	146.0&	1&	2&	2&	3&	3&	2&	117.4&	1&	3&	3&	3&	3&	3&	117.18\tabularnewline \hline
1&	1&	1&	3&	1&	1&	102.9&	1&	2&	2&	3&	2&	2&	82.68&	1&	3&	3&	3&	1&	3&	61.25\tabularnewline \hline
2&	3&	3&	3&	3&	3&	60.04&	2&	2&	3&	3&	2&	2&	52.84&	2&	3&	3&	3&	1&	3&	50.04\tabularnewline \hline
2&	1&	1&	3&	1&	1&	49.95&	2&	2&	3&	3&	1&	2&	49.60&	2&	1&	1&	3&	3&	1&	45.91\tabularnewline \hline
2&	3&	3&	3&	2&	3&	44.22&	2&	2&	3&	3&	3&	2&	35.86&	2&	1&	2&	3&	1&	1&	35.58\tabularnewline \hline
2&	2&	2&	3&	2&	2&	31.52&	2&	1&	1&	3&	2&	1&	31.38&	2&	1&	2&	3&	3&	1&	24.70\tabularnewline \hline

\end{tabular}
\end{table*}


\begin{table*}[hbt!]
   
\captionsetup{labelformat=empty}
\label{resulttable}
\centering
\begin{tabular}{|>{\centering}c |>{\centering}c |>{\centering}c |>{\centering}c|>{\centering}c|>{\centering}c |>{\centering}c|>{\centering}c |>{\centering}c|>{\centering}c|>{\centering}c |>{\centering}c |>{\centering}c |>{\centering}c|>{\centering}c|>{\centering}c |>{\centering}c |>{\centering}c |}

\hline
 \multicolumn{18}{|c|}{$HENK_{en,k}^{\acute{en},lev,t}$\;\;\;\; }\\
\hline

K&	en&	$\acute{en}$&	lev&	t&	Value&	kk&	en&	$\acute{en}$&	lev&	t&	Value&	kk&	en&	$\acute{en}$&	lev&	t&	Value\tabularnewline \hline
3&	1&	7&	1&	1&	6000&	2&	4&	4&	1&	1&	6000&	1&	2&	6&	1&	1&	60000\tabularnewline \hline
3&	1&	7&	2&	1&	400&	2&	4&	4&	2&	1&	400&	1&	2&	1&	2&	1&	400\tabularnewline \hline

\end{tabular}
\end{table*}




\begin{table*}[hbt!]
   
\captionsetup{labelformat=empty}
\label{resulttable}
\centering
\begin{tabular}{|>{\centering}c |>{\centering}c |>{\centering}c |>{\centering}c|>{\centering}c|>{\centering}c |>{\centering}c|>{\centering}c |>{\centering}c|>{\centering}c|>{\centering}c |>{\centering}c |>{\centering}c |>{\centering}c|>{\centering}c|>{\centering}c |>{\centering}c |>{\centering}c |>{\centering}c|>{\centering}c|>{\centering}c|}

\hline
 \multicolumn{21}{|c|}{$qlm_{p,l,m,v}^{t}, l \in L, m \in \acute{M}$\;\;\;\; (1000 barrel\;per\;day)}\\
\hline

p&	l&	m&	v&	t&	e&	value&	p&	l&	m&	v&	t&	e&	value&	p&	l&	m&	v&	t&	e&	value\tabularnewline \hline
1&	1&	1&	2&	3&	1&	180&	1&	2&	3&	1&	3&	2&	165.4&	1&	3&	2&	2&	3&	3&	162.\tabularnewline \hline
1&	3&	1&	2&	2&	3&	150.&	1&	1&	2&	2&	1&	1&	146.5&	1&	1&	3&	1&	2&	1&	127.\tabularnewline \hline
1&	2&	3&	1&	1&	2&	100&	&1	3&	1&	1&	1&	3&	73.25&	2&	3&	2&	1&	2&	3&	61.9\tabularnewline \hline
2&	3&	1&	1&	3&	3&	55.2&	1&	3&	2&	2&	1&	3&	53.41&	1&	2&	3&	1&	2&	2&	52.8\tabularnewline \hline
1&	1&	2&	2&	2&	1&	50.9&	2&	3&	1&	2&	1&	3&	50&	    2&	1&	3&	1&	1&	1&	49.9\tabularnewline \hline
1&	3&	2&	2&	2&	2&	49.1&	1&	2&	1&	2&	2&	2&	49.1&	2&	3&	1&	1&	2&	2&	44.8\tabularnewline \hline
2&	3&	3&	2&	3&	2&	35.8&	2&	2&	1&	1&	2&	1&	35.19&	1&	2&	2&	2&	3&	1&	34.5\tabularnewline \hline
2&	1&	2&	2&	3&	1&	31.7&	2&	1&	3&	1&	2&	1&	31.38&	2&	2&	2&	2&	3&	2&	28.5\tabularnewline \hline
1&	2&	1&	1&	1&	2&	26.7&	2&	2&	1&	1&	3&	1&	24.70&	1&	3&	1&	1&	3&	3&	20\tabularnewline \hline
1&	3&	3&	2&	3&	2&	14.5&	2&	1&	3&	2&	3&	1&	14.13&	2&	3&	2&	1&	3&	3&	9.65\tabularnewline \hline
2&	2&	3&	2&	2&	2&	8.61&	2&	3&	2&	1&	2&	2&	8.04&	1&	3&	2&	2&	3&	2&	2.86\tabularnewline \hline
2&	2&	2&	2&	1&	1&	0.35&	2&	3&	3&	2&	1&	3&	0.04&	&	&	&	&	&	&	\tabularnewline \hline

\end{tabular}
\end{table*}

\begin{table*}[hbt!]
   
\captionsetup{labelformat=empty}
\caption*{}
\label{resulttable}
\centering
\begin{tabular}{|>{\centering}c |>{\centering}c |>{\centering}c |>{\centering}c|>{\centering}c|>{\centering}c |>{\centering}c|>{\centering}c |>{\centering}c|>{\centering}c|>{\centering}c |>{\centering}c |>{\centering}c |>{\centering}c|>{\centering}c|>{\centering}c |>{\centering}c |>{\centering}c |>{\centering}c|>{\centering}c|>{\centering}c|>{\centering}c|>{\centering}c|>{\centering}c|>{\centering}c|>{\centering}c|>{\centering}c|>{\centering}c|}

\hline
 \multicolumn{24}{|c|}{$nlm_{p,l,m,v}^{lcv,t,e}, l \in L,m \in M$\;\;\;\; (Number)}\\
\hline

l&	m&	t&	lcv&	v&	p&	e&	value&	l&	m&	t&	lcv&	v&	p&	e&	value&	l&	m&	t&	lcv&	v&	p&	e&	value\tabularnewline \hline
2&	3&	3&	1&	1&	1&	2&	2&	3&	3&	3&	1&	2&	2&	2&	1&	3&	3&	3&	1&	2&	1&	2&	1\tabularnewline \hline
3&	3&	1&	1&	2&	2&	3&	1&	3&	2&	3&	1&	2&	1&	3&	1&	3&	2&	3&	1&	2&	1&	2&	1\tabularnewline \hline
3&	2&	3&	1&	1&	2&	3&	1&	3&	2&	2&	1&	2&	1&	2&	1&	3&	2&	2&	1&	1&	2&	3&	1\tabularnewline \hline
3&	2&	2&	1&	1&	2&	2&	1&	3&	2&	1&	1&	2&	1&	3&	1&	3&	2&	1&	1&	1&	2&	2&	1\tabularnewline \hline
3&	1&	3&	1&	1&	2&	3&	1&	3&	1&	3&	1&	1&	1&	3&	1&	3&	1&	3&	1&	1&	1&	3&	1\tabularnewline \hline
3&	1&	2&	1&	2&	1&	3&	1&	3&	1&	2&	1&	1&	2&	2&	1&	3&	1&	1&	1&	2&	2&	3&	1\tabularnewline \hline
3&	1&	1&	1&	2&	2&	3&	1&	3&	1&	1&	1&	1&	1&	3&	1&	2&	3&	2&	1&	2&	2&	2&	1\tabularnewline \hline
2&	3&	2&	1&	1&	1&	2&	1&	2&	3&	1&	1&	1&	1&	2&	1&	2&	2&	3&	1&	2&	2&	2&	1\tabularnewline \hline
2&	2&	3&	1&	2&	1&	1&	1&	2&	2&	1&	1&	2&	2&	1&	1&	2&	1&	3&	1&	1&	2&	1&	1\tabularnewline \hline
2&	1&	2&	1&	2&	1&	2&	1&	2&	1&	2&	1&	2&	1&	2&	1&	2&	1&	1&	1&	1&	1&	2&	1\tabularnewline \hline
1&	3&	3&	1&	2&	2&	1&	1&	1&	3&	2&	1&	1&	2&	1&	1&	1&	3&	2&	1&	1&	1&	1&	1\tabularnewline \hline
1&	3&	1&	1&	1&	2&	1&	1&	1&	2&	3&	1&	2&	2&	1&	1&	1&	2&	2&	1&	2&	1&	1&	1\tabularnewline \hline
1&	2&	1&	1&	2&	1&	1&	1&	1&	1&	3&	1&	2&	1&	1&	1&	&	&	&	&	&	&	&	\tabularnewline \hline

\end{tabular}

\end{table*}

\begin{table*}[hbt!]
   
\captionsetup{labelformat=empty}
\label{resulttable}
\centering
\begin{tabular}{|>{\centering}c |>{\centering}c |>{\centering}c |>{\centering}c|>{\centering}c|>{\centering}c |>{\centering}c|>{\centering}c |>{\centering}c|>{\centering}c|>{\centering}c |>{\centering}c |}

\hline
 \multicolumn{12}{|c|}{$Alab_{en}^{lev,t}$\;\;\;\; }\\
\hline

en&	lev&	t&	Value&	en&	lev&	t&	Value&	en&	lev&	t&	Value\tabularnewline \hline
3&	1&	3&	1.2e+5	&3&	1&	2	&80000&	1&	1&	2&	65158\tabularnewline \hline
2&	2&	3&	1.06e+5&	2&	2&	2	&70546&	8&	1&	3&	64800\tabularnewline \hline
8&	1&	2&	43200&	4&	1&	3&	41091&	3&	1&	1&	40000\tabularnewline \hline
5&	1&	3&	38532&	2&	1&	3&	36564&	1&	1&	1&	35579\tabularnewline \hline
2&	2&	1&	35473&	8&	2&	3&	30600&	6&	1&	3&	30333\tabularnewline \hline
1&	2&	3&	27032&	5&	1&	2&	25688&	3&	2&	3&	29712\tabularnewline \hline
4&	1&	2&	25394&	2&	1&	2&	22376&	8&	1&	1&	21600\tabularnewline \hline
8&	2&	2&	20400&	6&	1&	2&	20222&	3&	2&	2&	19808\tabularnewline \hline
1&	2&	2&	17888&	7&	1&	3&	16002&	4&	1&	1&	15697\tabularnewline \hline
2&	1&	1&	14188&	5&	1&	1&	12844&	4&	2&	3&	11381\tabularnewline \hline
9&	1&	3&	10701&	7&	1&	2&	10668&	9&	2&	3&	10599\tabularnewline \hline
8&	2&	1&	10200&	6&	1&	1&	10111&	3&	2&	1&	9904\tabularnewline \hline
5&	2&	3&	9633&	1&	2&	1&	9144&	6&	2&	3&	7581\tabularnewline \hline
4&	2&	2&	7454&	9&	1&	2&	7134&	9&	2&	2&	7066\tabularnewline \hline
5&	2&	2&	6422&	7&	1&	1&	5334&	6&	2&	2&	5054\tabularnewline \hline
7&	2&	3&	3990&	4&	2&	1&	3927&	9&	1&	1&	3567\tabularnewline \hline
9&	2&	1&	3533&	5&	2&	1&	3211&	7&	2&	2&	2660\tabularnewline \hline
6&	2&	1&	2527&	7&	2&	1&	1330&	&	&	&	\tabularnewline \hline

\end{tabular}
\end{table*}


\end{document}